\documentstyle[12pt]{article}
\newcommand{\sect}[1]{\section{#1}\setcounter{equation}{0}}
\newcommand{\subsect}[1]{\subsection{#1}}

\font\mbn=msbm10 scaled \magstep1
\font\mbs=msbm7 scaled \magstep1
\font\mbss=msbm5 scaled \magstep1
\newfam\mbff
\textfont\mbff=\mbn
\scriptfont\mbff=\mbs
\scriptscriptfont\mbff=\mbss\def\mbf{\fam\mbff}
\def\Re{{\mbf R}}
\def\Q{{\mbf Q}}
\def\Z{{\mbf Z}}
\def\Co{{\mbf C}}

\def\N{{\mbf N}}
\def\F{{\mbf F}}
\newtheorem{Th}{Theorem}[section]
\newtheorem{Lm}[Th]{Lemma}
\newtheorem{C}[Th]{Corollary}
\newtheorem{D}[Th]{Definition}
\newtheorem{Proposition}[Th]{Proposition}
\newtheorem{R}[Th]{Remark}
\newtheorem{E}[Th]{Example}
\author{Alexander Brudnyi\thanks{%Research supported in part by NSERC.
}\\
Department of Mathematics and Statistics\\
University of Calgary, Calgary\\
Canada}

\title{Formal Paths, Iterated Integrals and the Center Problem for Ordinary Differential Equations}

\date{} 

\begin{document} 
%==================================

%==================================

\maketitle

\begin{abstract}
{We continue the study of the center problem for the ordinary differential equation $v'=\sum_{i=1}^{\infty}a_{i}(x)v^{i+1}$ started in [Br1]-[Br5]. 
In this paper we present the highlights of the algebraic theory of centers.}
\end{abstract}
\tableofcontents

\sect{Introduction}
\quad
In this paper we describe an algebraic approach to the center problem for the ordinary differential equation
\begin{equation}\label{e1}
\frac{dv}{dx}=\sum_{i=1}^{\infty}a_{i}(x)v^{i+1},\ \ \ x\in I_{T}:=[0,T],
\end{equation}
with coefficients $a_{i}$ from the Banach space $L^{\infty}(I_{T})$ of bounded measurable complex-valued functions on $I_{T}$ equipped with the supremum norm. Condition
\begin{equation}\label{e2}
\sup_{x\in I_{T},\ i\in\N}\sqrt[i]{|a_{i}(x)|}<\infty
\end{equation}
guarantees that (\ref{e1}) has Lipschitz solutions on $I_{T}$ for all sufficiently small initial values.
By $X$ we denote the complex Fr\'{e}chet space of sequences $a=(a_{1},a_{2},\dots)$ satisfying (\ref{e2}).
We say that equation (\ref{e1}) determines a {\em center} if every solution $v$ of (\ref{e1}) with a sufficiently small initial value satisfies $v(T)=v(0)$. By ${\cal C}\subset X$ we denote the set of centers of (\ref{e1}). The center problem is: {\em given $a\in X$ to determine whether $a\in {\cal C}$}. It arises naturally in the framework of the geometric theory of ordinary differential equations created by Poincar\'{e}. In particular, there is a relation between the center problem for (\ref{e1}) and the classical Poincar\'{e} Center-Focus problem for planar polynomial vector fields 
\begin{equation}\label{e3}
\frac{dx}{dt}=-y+F(x,y),\ \ \ \frac{dy}{dt}=x+G(x,y),
\end{equation}
where $F$ and $G$ are polynomials of a given degree without constant and linear terms. This problem asks about conditions on $F$ and $G$ under which all trajectories of (\ref{e3}) situated in a small neighbourhood of $0\in\Re^{2}$ are closed. Passing to polar coordinates $(x,y)=(r\cos\phi, r\sin\phi)$ in (\ref{e3}) and expanding the right-hand side of the resulting equation as a series in $r$ (for $F$, $G$ with sufficiently small coefficients) we obtain an equation of the form (\ref{e1}) whose coefficients are trigonometric polynomials depending polynomially on the coefficients of (\ref{e3}). This reduces the Center-Focus Problem for (\ref{e3}) to the center problem for (\ref{e1}) with coefficients depending polynomially on a parameter. 

In this paper we continue the study of the center problem for equation (\ref{e1}) started in [Br1]-[Br5]. One of the basic objects of our approach is a
metrizable
topological group $G(X)$ defined by the coefficients of equations (\ref{e1}) (the, so-called, {\em group of paths} in $\Co^{\infty}$). Modulo the set of universal centers 
${\cal U}\subset {\cal C}$ of (\ref{e1}), described explicitly in [Br2], the set of centers forms a normal subgroup $\widehat{\cal C}\subset G(X)$. By
$G_{f}(X)$ and $\widehat{\cal C}_{f}$ we denote the groups of {\em formal paths} and of {\em formal centers}, respectively, i.e., the completions of 
$G(X)$ and $\widehat{\cal C}$ with respect to the metric on $G(X)$.
In this paper we study the algebraic properties of $G_{f}(X)$ and $\widehat{\cal C}_{f}$.
In particular, we describe Lie algebras of these groups and prove that
$G_{f}(X)$ is the semidirect product of a naturally defined normal subgroup of $\widehat {\cal C}_{f}$ and the 
subgroup $G_{f}(X^{2})$ of formal paths in $G_{f}(X)$ determined by coefficients of Abel differential equations, i.e., equations (\ref{e1}) with $a_{k}=0$ for all $k\geq 3$. Also, we show that $\widehat{\cal C}_{f}$ contains a dense subgroup of centers generated by certain piecewise linear paths in $\Co^{\infty}$. 

The paper is organized as follows.

Section 2 is devoted to the study of the group $G_{f}(X)$ of formal paths in $\Co^{\infty}$.

In section 2.1 we introduce a natural multiplication on the set $X$ similar to the product of paths in topology. Then we define the group of paths $G(X)$ as the quotient of $X$ by an equivalence relation determined in terms of iterated integrals on $X$. 
We equip $G(X)$ with a natural metric $d$ and define the group $G_{f}(X)$ of formal paths in $\Co^{\infty}$ as the completion of $G(X)$ with respect to $d$.

The equivalence class ${\cal U}\subset X$ corresponding to $1\in G(X)$ is called the set of universal centers of equation (\ref{e1}).
In section 2.2 we present some results of [Br2] on the characterization of elements from ${\cal U}$. 

Next, in section 2.3 we will show how to embed $G_{f}(X)$ in a group $G$ of invertible formal power series in $t$ whose coefficients belong to the associative algebra with unit $I$ of complex noncommutative polynomials in $I$ and free noncommutative variables $X_{i}$, $i\in\N$. 

Identifying $G_{f}(X)$ with its image in $G$ we describe in section 2.4 the Lie algebra ${\cal L}_{Lie}$ of $G_{f}(X)$ as the subset of Lie elements of the Lie algebra ${\cal L}_{G}$ of $G$. 

In section 2.5 we prove some structural theorems for $G(X)$ and $G_{f}(X)$. Namely, we describe the topological lower central series of these groups and their subgroups corresponding to closed paths in $\Co^{\infty}$ in terms of iterated integrals on $G_{f}(X)$. 

Finally, in section 2.6 we describe some natural subgroups of $G_{f}(X)$:
the groups $G_{f}(X^{n})$ generated by paths in $\Co^{n}$ and $G(X_{\F})$ determined over a field $\F\subset\Co$.

Section 3 is devoted to the study of the center problem for equation (\ref{e1}). 

In section 3.1 we gather some results from [Br1]-[Br3] on the explicit expression for the first return map of (\ref{e1}).  

Using this we define in section 3.2 the group $\widehat{\cal C}_{f}\subset G_{f}(X)$ of formal centers of equation (\ref{e1}).  

Then in section 3.3 we give an explicit description of the Lie algebra of $\widehat{\cal C}_{f}$ and show that $\widehat{\cal C}_{f}$ is the closure in $G_{f}(X)$ of the group of centers $\widehat{\cal C}$ of (\ref{e1}). We also prove that
$G_{f}(X)$ is the semidirect product of a naturally defined normal subgroup of $\widehat {\cal C}_{f}$ and the 
subgroup $G_{f}(X^{2})$ of formal paths in $\Co^{2}$ (i.e., determined by coefficients of Abel differential equations). At the end of this section we briefly discuss the center problem over a field $\F\subset\Co$. 

Finally, in section 3.4 we introduce the subgroup $PL\subset G_{f}(X)$ of piecewise linear paths in $\Co^{\infty}$. We give a characterization of centers belonging to this group and show that the set of such centers is dense in $\widehat{\cal C}_{f}$.
%===============
\sect{Group of Formal Paths}
\subsect{\hspace*{-1em}. Definition of the Group of Paths}

\quad
{\bf 2.1.1.} Let us consider $X$ as a semigroup with the operations given for $a=(a_{1},a_{2},\dots)$ and $b=(b_{1},b_{2},\dots)$ by
$$
a*b=(a_{1}*b_{1},a_{2}*b_{2},\dots)\in X\ \ \ {\rm and}\ \ \
a^{-1}=(a_{1}^{-1},a_{2}^{-1},\dots)\in X,
$$
where for $i\in\N$
$$
(a_{i}*b_{i})(x)=\left\{
\begin{array}{ccc}
2b_{i}(2x)&{\rm if}&0\leq x\leq T/2,\\
2a_{i}(2x-T)&{\rm if}&T/2<x\leq T
\end{array}
\right.
$$
and
$$
a_{i}^{-1}(x)=-a_{i}(T-x),\ \ \ 0\leq x\leq T.
$$

Let $\Co^{\infty}$ be the vector space of sequences of complex numbers $(c_{1},c_{2},\dots)$ equipped with the product topology.
For $a=(a_{1},a_{2},\dots)\in X$ by $\widetilde a= (\widetilde a_{1},
\widetilde a_{2},\dots): I_{T}\to\Co^{\infty}$, 
$\widetilde a_{k}(x):=\int_{0}^{x}a_{k}(t)\ \!dt$ for all $k\in\N$,
we denote a path in $\Co^{\infty}$ starting at 0.
The one-to-one map $a\mapsto\widetilde a$ sends the product $a*b$ to the product of paths $\widetilde a\circ\widetilde b$, that is, the path obtained by translating $\widetilde a$ so that its beginning meets the end of $\widetilde b$ and then forming the composite path. Similarly, $\widetilde {a^{-1}}$ is the path obtained by translating $\widetilde a$ so that its end meets 0 and then taking it with the opposite orientation. 
\\

{\bf 2.1.2.} For $a\in X$ let us consider the {\em basic iterated integrals}
\begin{equation}\label{e2.1}
I_{i_{1},\dots, i_{k}}(a):=\int\cdots\int_{0\leq s_{1}\leq\cdots\leq s_{k}\leq T}a_{i_{k}}(s_{k})\cdots a_{i_{1}}(s_{1})\ \!ds_{k}\cdots ds_{1}
\end{equation}
(for $k=0$ we assume that this equals $1$). By the Ree shuffle formula [R] the linear space generated by all such functions on $X$ is an algebra. The number $k$ in (\ref{e2.1}) is called the {\em order} of the iterated integral.
Also, the basic iterated integrals satisfy the following equations (see, e.g., [H, Propositions 2.9 and 2.12]):
\begin{equation}\label{e2.2}
I_{i_{1},\dots, i_{k}}(a*b)=I_{i_{1},\dots, i_{k}}(a)+\sum_{j=1}^{k-1}I_{i_{1},\dots,i_{j}}(a)\cdot I_{i_{j+1},\dots, i_{k}}(b)+I_{i_{1},\dots, i_{k}}(b).
\end{equation}
\begin{equation}\label{e2.3}
I_{i_{1},\dots, i_{k}}(a^{-1})=(-1)^{k}I_{i_{1},\dots, i_{k}}(a).
\end{equation}

For $a,b\in X$ we write $a\sim b$ if
all basic iterated integrals vanish at $a*b^{-1}$. Equations (\ref{e2.2})
and (\ref{e2.3}) imply that $a\sim b$ if and only if $I_{i_{1},\dots, i_{k}}(a)=I_{i_{1}\dots, i_{k}}(b)$ for all basic iterated integrals.
In particular, $\sim$ is an equivalence relation on $X$.
By $G(X)$ we denote the set of equivalence classes. Then $G(X)$ is a group
with the product induced by the product $*$ on $X$. By
$\pi:X\to G(X)$ we denote the map determined by the equivalence relation.
By the definition each iterated integral $I_{\cdot}$ is constant on fibres of
$\pi$ and therefore it determines a function $\widehat I_{\cdot}$ on $G(X)$ such that $I_{\cdot}=\widehat I_{\cdot}\circ\pi$. The functions $\widehat I_{\cdot}$ will be referred to as iterated integrals on $G(X)$. These functions separate the points on $G(X)$.

Next, we equip $G(X)$ with the weakest topology $\tau$ in which all basic iterated
integrals $\widehat I_{i_{1},\dots, i_{k}}$ are continuous. Then 
$(G(X),\tau)$ is a topological group. Moreover, $G(X)$ is metrizable: 
for $g, h\in G(X)$ the formula
\begin{equation}\label{e2.4}
d(g,h):=\sum_{n=1}^{\infty}\frac{1}{4^{n}}\cdot\left(\ \!\sum_{i_{1}+\cdots +i_{k}=n}\frac{|\widehat I_{i_{1},\dots, i_{k}}(g)-\widehat I_{i_{1},\dots, i_{k}}(h)|}{1+|\widehat I_{i_{1},\dots, i_{k}}(g)-\widehat I_{i_{1},\dots, i_{k}}(h)|}\right)
\end{equation}
determines a metric on $G(X)$ compatible with topology $\tau$ (see [Br3, Theorem 2.4]). We mention also that $G(X)$ is contractible, residually torsion free nilpotent (i.e., finite dimensional unipotent representations of $G(X)$ separate the points on $G(X)$) and is the union of an increasing sequence of compact subsets (see [Br3]).

By $G_{f}(X)$ we denote the completion of $G(X)$ with respect to the metric
$d$. Then $G_{f}(X)$ is a topological group which will be called the {\em group of formal paths} in $\Co^{\infty}$.
%==========
\subsect{\hspace*{-1em}. Structure of the Set of Universal Centers}

\quad By ${\cal U}\subset X$ we denote the set of elements $a\in X$ such that
all basic iterated integrals vanish at $a$. According to (\ref{e2.2}), 
${\cal U}$ is a sub-semigroup of $X$. It was shown in [Br2] that ${\cal U}\subset {\cal C}$, the set of centers of equation (\ref{e1}). We call ${\cal U}$ 
the {\em set of universal centers} of (\ref{e1}). In this section we formulate
some results on the characterization of elements from ${\cal U}$ established in [Br2].

Let $X^{k}:=\{a=(a_{1},a_{2},\dots)\in X\ :\ a_{j}=0\ {\rm for\ all}\
j>k\}$. By $p_{k}: X\to X^{k}$, $(a_{1},a_{2},\dots)\mapsto (a_{1},\dots, a_{k},0,0,\dots)$, we denote the natural projection. Clearly $a\in {\cal U}$ if
and only if $p_{k}(a)\in {\cal U}$ for all $k\in\N$. Therefore it suffices to characterize elements from the sets ${\cal U}^{\ \!\!k}:={\cal U}\cap X^{k}$.

For $a=(a_{1},\dots, a_{k},0,\dots)\in X^{k}$ consider the Lipschitz curve
$A_{k}:I_{T}\to\Co^{k}$ determined by the formula
\begin{equation}\label{e2.5}
A_{k}(x):=\left(\int_{0}^{x}a_{1}(t)\ \!dt,\dots,\int_{0}^{x}a_{k}(t)\ \!dt\right),\ \ \
x\in I_{T}.
\end{equation}
We set $\Gamma_{k}:=A_{k}(I_{T})$.

Next, we require
\begin{D}\label{d2.1}
The polynomially convex hull $\widehat K$ of a compact set $K\subset\Co^{k}$ is the set of points $z\in\Co^{k}$ such that if $p$ is any holomorphic polynomial in $k$ variables
$$
|p(z)|\leq\max_{x\in K}|p(x)|.
$$
\end{D}

It is well known (see, e.g., [AW]) that $\widehat K$ is compact, and if $K$ is connected then $\widehat K$ is connected. The following basic result was proved in [Br2, Theorem 1.10].
\begin{Th}\label{te2.2}
Suppose that $a\in {\cal U}^{\ \!\!k}$. Then for any domain $U\subset\Co^{k}$
containing $\widehat\Gamma_{k}$ the path $A_{k}:I_{T}\to U$ is closed and represents the unit element of the fundamental group $\pi_{1}(U)$ of $U$.
\end{Th}

Since $A_{k}$ is Lipschitz, $\Gamma_{k}$ is of a finite linear measure. Then according to the result of Alexander [A], $\widehat{\Gamma}_{k}\setminus
\Gamma_{k}$ is a (possibly empty) pure one-dimensional complex analytic subset of $\Co^{k}\setminus\Gamma_{k}$. In particular, since the covering dimension
of $\Gamma_{k}$ is $1$, the covering dimension of $\widehat{\Gamma}_{k}$ is 2. Therefore according to the Freudenthal expansion theorem [F], $\widehat\Gamma_{k}$ can be
presented as an inverse limit of a sequence
of compact polyhedra $\{Q_{kj}\}_{j\in\N}$ with $dim\ \!Q_{kj}\leq 2$.
Let $\pi_{kj}:\widehat{\Gamma}_{k}\to Q_{kj}$ be the
continuous projections generated by the inverse limit construction. 
It is easy to check that Theorem \ref{te2.2} is equivalent to the following
statement:  {\em if $a\in {\cal U}^{\ \!\!k}$, then for all $j$ 
the continuous
paths $\pi_{kj}\circ A_{k}:I_{T}\to Q_{kj}$ are closed 
and represent unit elements of $\pi_{1}(Q_{kj})$}.

Let us formulate two corollaries of Theorem \ref{te2.2}. In the first one we use the notion of a {\em bordered Riemann surface}. This is a compact connected set which consists of a (possibly singular) one-dimensional complex analytic space with a $C^{2}$-boundary.
\begin{C}\label{c2.3}
Suppose that $a\in X^{k}$ is such that the corresponding set $\widehat\Gamma_{k}$ belongs to a bordered Riemann surface $S\subset\Co^{k}$. Then $a\in {\cal U}^{\ \!\!k}$ if and only if the path
$A_{k}:I_{T}\to S$ is closed and represents the unit element of the fundamental group $\pi_{1}(S)$ of $S$.
\end{C}

The proof of Corollary \ref{c2.3} repeats literally the arguments of the proofs of [Br2, Corollary 1.17] and [Br5, Corollary 3.7]. Using the covering homotopy theorem one obtains the following reformulation of the above result.

Let $\pi:\widetilde S\to S$ be the universal covering of $S$. {\em Under the hypotheses of Corollary \ref{c2.3}, $a\in {\cal U}^{\ \!\!k}$ if and only if there is a closed Lipschitz path $\widetilde A_{k}:I_{T}\to\widetilde S$ such that $A_{k}=\pi\circ\widetilde A_{k}$}. 
\begin{E}\label{ex2.4}
{\rm 1) Suppose that $\Gamma_{k}\subset\Co^{k}$ is the image of the unit circle under a holomorphic embedding of its neighbourhood.
Then by the result of Wermer [W], $\widehat{\Gamma}_{k}$ is a bordered Riemann surface.\\
(2) The assumptions of Corollary \ref{c2.3} are also fulfilled if $\Gamma_{k}$ belongs to a one-dimensional complex analytic subset of a domain $U\subset\Co^{k}$ such that $U=\cup_{j}K_{j}$ with $K_{j}\subset\subset K_{j+1}$ and $\widehat{K}_{j}=K_{j}$ for all $j\in\N$, cf. [Br2, Corollary 1.17]. In particular, this is valid for a convex $U$.}
\end{E}

To formulate the second corollary of Theorem \ref{te2.2} we recall the definition of a Lipschitz triangulable curve, see, e.g., [Br2]. 
\begin{D}\label{d2.5}
A compact curve $C\subset\Re^{N}$ is called Lipschitz triangulable if
\begin{enumerate}
\item
$C=\cup_{j=1}^{s}C_{j}$ and for $i\neq j$ the intersection $C_{i}\cap C_{j}$ consists of at most one point;
\item
There are Lipschitz embeddings $f_{j}:[0,1]\to\Re^{N}$ such that $f_{j}([0,1])=C_{j}$;
\item
The inverse maps $f_{j}^{-1}:C_{j}\to\Re$ are locally Lipschitz on $C_{j}\setminus \{f_{j}(0)\cup f_{j}(1)\}$.
\end{enumerate}
\end{D}

The following corollary extends one of the main results of Chen [C3].
\begin{C}\label{c2.6}
Suppose that $a\in X^{k}$ is such that
$\Gamma_{k}$ is a Lipschitz triangulable curve and $\widehat{\Gamma}_{k}=\Gamma_{k}$. Then $a\in  {\cal U}^{\ \!\!k}$
if and only if the path $A_{k}:I_{T}\to\Gamma_{k}$ is closed and represents
the unit element of $\pi_{1}(\Gamma_{k})$.
\end{C}

Using the covering homotopy theorem one reformulates this corollary as follows:

{\em Under the hypotheses of Corollary \ref{c2.6}, $a\in  {\cal U}^{\ \!\!k}$ if and only if there are a Lipschitz triangulable curve ${\cal T}$ homeomorphic to a finite tree, a locally bi-Lipschitz map $\pi:{\cal T}\to\Gamma_{k}$ and a closed Lipschitz path
$\widetilde A_{k}:I_{T}\to {\cal T}$ such that $A_{k}=\pi\circ\widetilde A_{k}$}.

\begin{E}\label{ex2.7}
{\rm (1) The condition $\widehat{\Gamma}_{k}=\Gamma_{k}$ is fulfilled if, e.g.,  $\Gamma_{k}$ belongs to a compact set $K$ in a $C^{1}$-manifold $M$ with no complex tangents such that $\widehat K=K$ (for the proof see, e.g., [AW, Theorem 17.1]). For instance, one can take as such $K$ any compact subset of $M=\Re^{k}$.\\
(2) $\Gamma_{k}$ is a Lipschitz triangulable curve if, e.g., $A_{k}:I_{T}\to\Co^{k}$ is nonconstant analytic.}
\end{E}

Corollaries \ref{c2.3} and \ref{c2.6} reveal a connection of the center problem for equation (\ref{e1}) with the so-called {\em composition condition} whose role and importance was studied in [AL], [BFY1], [BFY2], [Y] for the special case of Abel differential equations.
%============
\subsect{\hspace*{-1em}. Representation of Paths by Noncommutative Formal Power Series}

\quad {\bf 2.3.1.}
Let $\Co\!\left\langle X_{1},X_{2},\dots\right\rangle$ be the associative algebra with unit $I$ of complex noncommutative polynomials in $I$ and free noncommutative variables $X_{1}, X_{2},\dots$ (i.e., there are no nontrivial relations between these variables). By $\Co\!\left\langle X_{1},X_{2},\dots\right\rangle\![[t]]$ we denote the associative algebra of formal power series in $t$ with coefficients from
$\Co\!\left\langle X_{1},X_{2},\dots\right\rangle$. Also, by ${\cal A}\subset \Co\!\left\langle X_{1},X_{2},\dots\right\rangle\![[t]]$ we denote the subalgebra of series $f$ of the form
\begin{equation}\label{e2.6}
f=c_{0}I+\sum_{n=1}^{\infty}\left(\ \!\sum_{i_{1}+\dots +i_{k}=n}c_{i_{1},\dots, i_{k}}X_{i_{1}}\cdots X_{i_{k}}
\right)t^{n}
\end{equation}
with $c_{0},c_{i_{1},\dots, i_{k}}\in\Co$ for all $i_{1},\dots, i_{k}, k\in\N$.
\begin{R}\label{re2.8}
{\rm Let $S\subset\Co\!\left\langle X_{1},X_{2},\dots\right\rangle$ be the multiplicative semigroup generated by  $I,X_{1},X_{2},\dots$. Consider a grading function $w:S\to\Z_{+}$ determined
by the conditions
$$
w(I)=0,\ w(X_{i})=i,\ \ \ i\in\N,\ \ \ {\rm and}\ \ \ w(x\cdot y):=w(x)+w(y)\ \ \
{\rm for\ all}\ \ \ x,y\in S.
$$
This splits $S$ in a disjoint union $S=\sqcup_{n=0}^{\infty}S_{n}$, where
$S_{n}=\{s\in S\ :\ w(s)=n\}$. Now, each ${f\in\cal A}$ is written as 
\begin{equation}\label{e2.7}
f=\sum_{n=0}^{\infty}f_{n} t^{n}\ \ \ {\rm where}\ \ \ f_{n}\in V_{n}:=span_{\Co}(S_{n}),\  n\in\Z_{+}.
\end{equation}
}
\end{R}

We equip ${\cal A}$ with the weakest topology in which all coefficients 
$c_{i_{1},\dots, i_{k}}$ in (\ref{e2.6}) considered as functions in $f\in {\cal A}$ are continuous.
Since the set of these functions is countable, ${\cal A}$ is metrizable,
cf. (\ref{e2.4}). 
Moreover, if $d$ is a metric on ${\cal A}$ compatible with the topology, then $({\cal A}, d)$ is a complete metric space. Also, by the definition the multiplication $\cdot :{\cal A}\times {\cal A}\to {\cal A}$ is continuous in this topology. 
\begin{R}\label{r2.9.}
{\rm A sequence $\{f_{k}=\sum_{n=0}^{\infty}f_{nk}t^{n}\ :\ f_{nk}\in V_{n},\ n\in\Z_{+}\}_{k\in\N}$
converges to $f=\sum_{n=0}^{\infty}f_{n}t^{n}$, $f_{n}\in V_{n}$, $n\in\Z_{+}$, in the topology of ${\cal A}$ if and only if each 
$\{f_{kn}\}_{k\in\N}$ converges to $f_{n}$ in $V_{n}$ naturally identified with the hermitian space $\Co^{d(n)}$ where $d(n):=\# S_{n}$.}
\end{R} 

By $G\subset{\cal A}$ we denote the closed subset of elements $f$ of form (\ref{e2.6}) with $c_{0}=1$. Then $(G,\cdot)$ is a topological group. Its Lie algebra ${\cal L}_{G}\subset {\cal A}$ consists of elements $f$ of form (\ref{e2.6}) with $c_{0}=0$. (For $f,g\in {\cal L}_{G}$ their product is defined by the formula $[f,g]:=f\cdot g-g\cdot f$.) Also, the map $\exp:{\cal L}_{G}\to G$,
$\exp(f):=e^{f}=\sum_{n=0}^{\infty}\frac{f^{n}}{n!}$, is a homeomorphism.
\\

{\bf 2.3.2.} For an element
$a=(a_{1},a_{2},\dots)\in X$ let us consider the equation 
\begin{equation}\label{e2.8}
F'(x)=\left(\ \!\sum_{i=1}^{\infty}a_{i}(x)\ \!t^{i}X_{i}\right)F(x),
\ \ \ \ \ x\in I_{T}.
\end{equation}
This can be solved by Picard iteration to obtain a solution $F_{a}: I_{T}\to G$, $F_{a}(0)=I$, whose coefficients in expansion in $X_{1}, X_{2},\dots$ and $t$ are Lipschitz functions on $I_{T}$.  We set 
\begin{equation}\label{e2.9}
E(a):=F_{a}(T),\ \ \ a\in X.
\end{equation}
By the definition, see section 2.1.1, we have, cf. [C1, Theorem 6.1],
\begin{equation}\label{e2.10}
E(a*b)=E(a)\cdot E(b),\ \ \ a,b\in X.
\end{equation}
Also, an explicit calculation leads to the formula
\begin{equation}\label{e2.11}
E(a)=I+\sum_{n=1}^{\infty}\left(\ \!\sum_{i_{1}+\cdots +i_{k}=n}I_{i_{1},\dots, i_{k}}(a)X_{i_{1}}\cdots X_{i_{k}}\right)t^{n}.
\end{equation}

The last formula implies that the kernel of the homomorphism $E:X\to G$ is the set of universal centers ${\cal U}$. In particular, there is a homomorphism $\widehat E:G(X)\to G$ such that $E=\widehat E\circ\pi$, that is,
\begin{equation}\label{e2.12}
\widehat E(g)=I+\sum_{n=1}^{\infty}\left(\ \!\sum_{i_{1}+\cdots +i_{k}=n}\widehat I_{i_{1},\dots, i_{k}}(g)X_{i_{1}}\cdots X_{i_{k}}\right)t^{n},\ \ \ g\in G(X). 
\end{equation}

Formula (\ref{e2.12}) shows that $\widehat E:G(X)\to G$ is a continuous embedding. Moreover, one can determine a metric $d_{1}$ on ${\cal A}$ compatible with topology such that $\widehat E:
(G(X),d)\to (G,d_{1})$ is an isometric embedding, cf. (\ref{e2.4}). Therefore $\widehat E$ is naturally extended to a continuous embedding $G_{f}(X)\to G$
(denoted also by $\widehat E$). By the definition, $\widehat E:G_{f}(X)\to G$ is an injective homomorphism of topological groups and $\widehat E(G_{f}(X))$ is
the closure of $\widehat E(G(X))$ in the topology of $G$. 

In what follows we identify $G(X)$ and $G_{f}(X)$ with their images under $\widehat E$. 
%=============
\subsect{\hspace*{-1em}. Lie Algebra of the Group of Formal Paths}

\quad {\bf 2.4.1.} Recall that each element $g\in {\cal L}_{G}$ 
can be written as
$g=\sum_{n=1}^{\infty}g_{n}t^{n}$, $g_{n}\in V_{n}$, $n\in\N$. We say that such $g$ is a {\em Lie element} if each $g_{n}$ belongs to the free Lie algebra generated by $X_{1},\dots, X_{n}$. In this case each $g_{n}$ has
the form
\begin{equation}\label{e2.13}
g_{n}=\sum_{i_{1}+\dots +i_{k}=n}c_{i_{1},\dots, i_{k}}[X_{i_{1}},[X_{i_{2}},[\ \cdots , [X_{i_{k-1}},X_{i_{k}}]\cdots \ ]]].
\end{equation}
with all $c_{i_{1},\dots, i_{k}}\in\Co$. (Here the term with $i_{k}=n$ is $c_{n}X_{n}$.) 

Let $L_{n}\subset V_{n}$ be the subspace of elements $g_{n}$ of form (\ref{e2.13}). It follows from 
[M-KO, Theorem 3.2] that  
\begin{equation}\label{e2.14}
dim_{\Co}L_{n}=\frac{1}{n}\sum_{d|n}(2^{\ \!n/d}-1)\cdot\mu(d)
\end{equation}
where the sum is taken over all numbers $d\in\N$ that divide $n$, and $\mu:\N\to\{-1,0,1\}$ is the M\"{o}bius function defined as follows.
If $d$ has a prime factorization
$$
d=p_{1}^{n_{1}}p_{2}^{n_{2}}\cdots p_{q}^{n_{q}},\ \ \ n_{i}>0,
$$
then
$$
\mu(d)=
\left\{
\begin{array}{cc}
1&{\rm for}\ d=1\\
\\
(-1)^{q}&\quad {\rm if\ all}\ n_{i}=1\\
\\
0&\ {\rm otherwise}.
\end{array}
\right.
$$

By ${\cal L}_{Lie}$ we denote the subset of Lie elements of ${\cal L}_{G}$. Then ${\cal L}_{Lie}$ is a closed (in the topology of ${\cal A}$) Lie subalgebra of ${\cal L}_{G}$.
\begin{Th}\label{te2.10}
The exponential map $\exp: {\cal L}_{G}\to G$ maps ${\cal L}_{Lie}$ homeomorphically onto $G_{f}(X)$.
\end{Th}

Thus ${\cal L}_{Lie}$ can be considered as the Lie algebra of $G_{f}(X)$.
\\

{\bf 2.4.2. Proof of Theorem \ref{te2.10}.} Let $\log :G\to {\cal L}_{G}$,
$\log(f)=-\sum_{n=1}^{\infty}\frac{(I-f)^{n}}{n}$, be the logarithmic map. By the definition it is continuous and inverse to the exponential map $\exp$.
From [C2, Theorem 4.2] and [R] follow that $\log$ maps $G(X)$ into
${\cal L}_{Lie}$. Since ${\cal L}_{Lie}$ is a closed subspace of
${\cal L}_{G}$ and $\log$ is continuous, it maps $G_{f}(X)$ into ${\cal L}_{Lie}$, as well. In particular, $G_{f}(X)\subset\exp({\cal L}_{Lie})$.
Let us prove the converse implication.

Let $J\subset {\cal A}$ be the two-sided ideal of elements $f$ of form (\ref{e2.6}) with $c_{0}=0$. By $J^{l}$ we denote the $l$th power of $J$.
Let $q_{l}:{\cal A}\to {\cal A}/J^{l}=:A_{l}$ be the quotient homomorphism.
We set $\overline{X}_{s}=q_{l}(X_{s}\cdot t^{s})$, $1\leq s\leq l-1$. Then  for
$f\in {\cal A}$ of form (\ref{e2.6})
we have
\begin{equation}\label{e2.15}
q_{l}(f):=c_{0}I+\sum_{n=1}^{l-1}\left(\ \!\sum_{i_{1}+\cdots +i_{k}=n}
c_{i_{1},\dots, i_{k}}\overline{X}_{i_{1}}\cdots\overline{X}_{i_{k}}\right).
\end{equation}
(Here $I$ is the unit of $A_{l}$.) We naturally identify $A_{l}$ with the
hermitian space $\Co^{n(k)}$, $n(k)=dim_{\Co}A_{l}$, so that
$q_{l}$ is a continuous map. Then $G_{l}:=q_{l}(G)\subset A_{l}$ is a complex nilpotent Lie group. 

Further, let $X_{rect}\subset X$ be the sub-semigroup of rectangular paths, i.e., elements $a\in X$ whose first integrals 
$\widetilde a: I_{T}\to\Co^{\infty}$ are paths consisting of segments each going in the direction of some particular coordinate. By $G(X_{rect})\subset G(X)$ we denote the subgroup generated by $X_{rect}$. Identifying $G(X)$ with its image in $G$ by $\widehat E$ we obtain from (\ref{e2.11}) that $G(X_{rect})$ is a subgroup of $G$ generated by elements $e^{c_{n}X_{n}t^{n}}$, $c_{n}\in\Co$, $n\in\N$. (In particular, $G(X_{rect})$ is isomorphic to the free product of countably many copies of $\Co$.) 
\begin{Proposition}\label{pr2.11}
The images of $G(X_{rect})$, $G(X)$ and $G_{f}(X)$ in $G_{l}$ coincide and form a complex Lie subgroup of $G_{l}$.
\end{Proposition}
{\bf Proof.}
We set for brevity
\begin{equation}\label{e2.16}
\overline{Q}_{l}:=q_{l}(G(X_{rect})),\ \ \ Q_{l}:=q_{l}(G(X)),\ \ \
\widetilde Q_{l}:=q_{l}(G_{f}(X)).
\end{equation}

For any $c\in\Co$, $a=(a_{1},a_{2},\dots)\in X$ by $ca$ we denote the element
$(ca_{1},c^{2}a_{2},\dots)\in X$. Suppose that $S\subset X$ is a sub-semigroup such that for each $s\in S$ and any $c\in\Co$ elements $s^{-1}, cs$ belong to $S$. Let $S_{l}:=(q_{l}\circ\widehat E\circ\pi)(S)$ be the image of $S$ in $G_{l}$. By the definition $S_{l}$ is a subgroup of $G_{l}$.
We will use the following 
\begin{Lm}\label{le2.12}
$S_{l}$ is a complex Lie subgroup of $G_{l}$.
\end{Lm}
{\bf Proof.}
Let $g\in G(X)\subset {\cal A}$ be the image of an element $a\in X$. For any $c\in\Co$ by $cg\in G(X)$ we denote the image of $ca\in X$. 
If $g=\sum_{n=0}^{\infty}g_{n}t^{n}$, $g_{n}\in V_{n}$, see (\ref{e2.7}), then
\begin{equation}\label{eq2.17}
q_{l}(g)=I+\sum_{n=1}^{l-1}q_{n}(g_{n}t^{n})\ \ \ {\rm and}\ \ \
q_{l}(cg):=I+\sum_{n=1}^{l-1}c^{n}\cdot q_{n}(g_{n}t^{n}).
\end{equation}
We will naturally identify $G_{l}$ with $\Co^{N}$ where 
$N:=dim_{\Co}\ \!G_{l}$. 

Let $K=\{g_{1},\dots, g_{k}\}\subset \pi(S)\subset G(X)$ be a finite set.
We define a holomorphic polynomial map $F_{K}:\Co^{k}\to\Co^{N}$ by  the formula
\begin{equation}\label{eq2.18}
F_{K}(z_{1},\dots, z_{k}):=q_{l}((z_{1}g_{1})\cdots (z_{k}g_{k})).
\end{equation}
Let $Z_{K}$ be the Zariski closure of $F_{K}(\Co^{k})$ in $\Co^{N}$. Then $Z_{K}$ is an irreducible complex algebraic subvariety of $\Co^{N}$ and $F_{K}(\Co^{k})$ contains an open dense subset of $Z_{K}$ (see, e.g., the book of Mumford [M] for the basic facts of Algebraic Geometry). By the definition we have 
\begin{equation}\label{eq2.19}
F_{K_{1}}(\Co^{\#K_{1}})\subset F_{K_{2}}(\Co^{\#K_{2}})\ \ \ {\rm and}\ \ \
 Z_{K_{1}}\subset Z_{K_{2}}\ \ \ {\rm for}\ \ \ K_{1}\subset K_{2}. 
\end{equation}
Let $Z:=\cup_{K}Z_{K}$ where $K$ runs over all finite subsets of $G_{l}$. Since $dim_{\Co}\ \!Z_{K}\leq N$ and all $Z_{K}$ are irreducible, (\ref{eq2.19}) implies that there is a finite set $K\subset\pi(S)$ such that $Z_{K}=Z$. Observe that by the definition the group
$S_{l}$ is dense in $Z=Z_{K}$. Hence $Z$ is the Zariski closure of $S_{l}$. In particular,
$Z$ is a complex Lie subgroup of $G_{l}$. Also, from the identity $Z=Z_{K}$ it follows that $S_{l}$ contains an open dense 
subset of $Z$. Since the topologies of 
the groups $S_{l}$ and $Z$ are induced from that of $G_{l}$, the latter implies that $S_{l}=Z$ completing the proof of the lemma.\ \ \ \ \ $\Box$

Continuing the proof of the proposition we
choose as the $S$ in Lemma \ref{le2.12} semigroups $X_{rect}$ and $X$. Then we conclude that
$\overline{Q}_{l}\subset Q_{l}$ are complex Lie subgroups of $G_{l}$ (in particular, they are closed subsets of $G_{l}$). 
Since $Q_{l}$ is dense in $\widetilde Q_{l}$, see section 2.1.2, the latter implies that $Q_{l}=\widetilde Q_{l}$.

Let ${\cal L}_{G_{l}}$ be the Lie algebra of $G_{l}$. By the definition it consists of all elements of form (\ref{e2.15}) with $c_{0}=0$.
Let $\exp :{\cal L}_{G_{l}}\to G_{l}$, $\exp(f):=e^{f}$, be the corresponding exponential map. Clearly we have the following commutative diagram:
\begin{equation}\label{e2.20}
\begin{array}{ccc}
\ \ \ {\cal L}_{G}&\stackrel{\exp}{\longrightarrow}&\!\!\!G\\
q_{l}\downarrow& &\ \downarrow q_{l}\\
\ \ \ \ \ {\cal L}_{G_{l}}&\stackrel{\exp}{\longrightarrow}& G_{l}.
\end{array}
\end{equation}
This implies that for each $f\in Q_{l}$ the element $\log(f):=-\sum_{n=1}^{l-1}\frac{(I-f)^{n}}{n}$ belongs to the free nilpotent Lie algebra
${\cal L}_{Lie}^{l}\subset {\cal L}_{G_{l}}$ generated by elements
$\overline{X}_{1},\dots,\overline{X}_{l-1}$, i.e., each element $g$ of this algebra has the form
\begin{equation}\label{e2.21}
g=\sum_{n=1}^{l-1}\left(\ \!\sum_{i_{1}+\dots +i_{k}=n}c_{i_{1},\dots, i_{k}}[\overline{X}_{i_{1}},[\overline{X}_{i_{2}},[\ \cdots , [\overline{X}_{i_{k-1}},\overline{X}_{i_{k}}]\cdots \ ]]]\right)
\end{equation}
with all $c_{i_{1},\dots, i_{k}}\in\Co$, $i_{1},\dots, i_{k}\in\{1,\dots, l-1\}$. Thus the Lie algebra ${\cal L}_{Q_{l}}:=\log(Q_{l})$ of $Q_{l}$ is a subset of ${\cal L}_{Lie}^{l}$. Next, the Lie algebra
${\cal L}_{\overline{Q}_{l}}:=\log(\overline{Q}_{l})$ of $\overline{Q}_{l}$ is a subset of ${\cal L}_{Lie}^{l}$, as well. Moreover, since $\overline{Q}_{l}$ is generated by elements $e^{c_{s}\overline{X}_{s}}$, $c_{s}\in\Co$, $1\leq s\leq l-1$, its Lie algebra is generated by elements $\overline{X}_{s}$, $1\leq s\leq l-1$, and therefore coincides with ${\cal L}_{Lie}^{l}$. This implies that ${\cal L}_{\overline{Q}_{l}}={\cal L}_{Q_{l}}={\cal L}_{Lie}^{l}$, and $Q_{l}=\overline{Q}_{l}$.

The proof of the proposition is complete.\ \ \ \ \ $\Box$

Let us finish the proof of the theorem. Let $g\in {\cal L}_{Lie}$. Consider the elements $g_{l}:=q_{l}(g)\in {\cal L}_{Lie}^{l}$, $l\geq 2$.
According to Proposition \ref{pr2.11} there are elements $f_{l}\in G(X_{rect})$ such that $q_{l}(f_{l})=e^{g_{l}}$. Now, for $m>l$ we have
$$
q_{l}(f_{m}\cdot f_{l}^{-1})=q_{l}(f_{m})\cdot (q_{l}(f_{l}))^{-1}=
q_{l}(q_{m}(f_{m}))\cdot e^{-g_{l}}=q_{l}(e^{g_{m}})\cdot e^{-g_{l}}=
e^{g_{l}}\cdot e^{-g_{l}}=I.
$$
Thus $f_{m}-f_{l}\in J^{l}$.
By the definition of the topology of $G$ this implies that the sequence
$\{f_{l}\}_{l\geq 2}$ converges in $G$ to an element $f\in G_{f}(X)$, so that
$q_{l}(f)=e^{g_{l}}$, $l\geq 2$. Taking here the limit as $l$ tends to $\infty$ we get $f=e^{g}$. This shows that $\exp({\cal L}_{Lie})\subset G_{f}(X)$ and completes the proof of the theorem.\ \ \ \ \ $\Box$
\begin{R}\label{r2.13}
{\rm We also established in the proof that $G(X_{rect})$ is a dense subgroup of $G_{f}(X)$.}
\end{R}

{\bf 2.4.3. Shuffles.}
\begin{D}\label{d2.14}
A permutation $\sigma$ of $\{1,2,\dots,r+s\}$ is a shuffle of type 
$(r,s)$ if
$$
\sigma^{-1}(1)<\sigma^{-1}(2)<\cdots <\sigma^{-1}(r)
$$
and
$$
\sigma^{-1}(r+1)<\sigma^{-1}(r+2)<\cdots <\sigma^{-1}(r+s).
$$
\end{D}
(The term ``shuffle`` is used because such permutations arise in riffle shuffling a deck of $r+s$ cards cut into one pile of $r$ cards and a second pile of $s$ cards.)

The following result is a corollary of Theorem \ref{te2.10}.
\begin{Th}\label{te2.15}
An element
$$
f=I+\sum_{n=1}^{\infty}\left(\ \!\sum_{i_{1}+\dots +i_{k}=n}c_{i_{1},\dots, i_{k}}X_{i_{1}}\cdots X_{i_{k}}\right) t^{n}\in G
$$
belongs to $G_{f}(X)$ if and only if its coefficients satisfy the system of Ree shuffle relations:
\begin{equation}\label{e2.22}
c_{i_{1},\dots , i_{r}}\cdot c_{i_{r+1},\dots, i_{r+s}}=\sum_{\sigma}c_{i_{\sigma(1)},\dots , i_{\sigma(r+s)}} ,\ \ \
i_{1},\dots, i_{r+s}\in\N,
\end{equation}
where the sum is taken over the set of shuffles of type $(r,s)$.
\end{Th}
{\bf Proof.} According to the main result of Ree [R], $\log(f)\in {\cal L}_{Lie}$ for $f\in G$ if and only if the coefficients of $f$ satisfy equations (\ref{e2.22}). This and Theorem \ref{te2.10} imply the required result.\ \ \ \ \ $\Box$
%=============
\subsect{\hspace*{-1em}. Topological Lower Central Series of Some Groups of Paths}
%Structure of the Group of Formal Paths}

\quad {\bf 2.5.1.} In the section we describe the topological lower central series of groups $G(X)$ and $G_{f}(X)$. 

Let $G$ be a topological group. We set $G_{n}:=[G,G_{n-1}]$ and
$G_{1}=[G,G]$ (the commutator subgroup of $G$). For $H\subset G$ by $\overline{H}\subset G$ we denote the closure of $H$.
\begin{D}\label{d2.16}
The sequence 
$$
G\supset\overline{G}_{1}\supset\overline{G}_{2}\supset\cdots
$$
is called the topological lower central series of $G$.
\end{D}

Next, consider the family $\{\widehat I_{i_{1},\dots, i_{k}}\}$ of all basic integrals on $G(X)$, see section 2.1.2. By the definition each function of this family admits a continuous extension to $G_{f}(X)$. We retain the same symbols for the extended functions and call them the {\em basic iterated integrals} on $G_{f}(X)$. Observe that if $c_{i_{1},\dots, i_{k}}:{\cal A}\to\Co$ is the function whose value at $f\in {\cal A}$ is the coefficient
corresponding to the monomial $X_{i_{1}}\cdots X_{i_{k}}$ in the series expansion (\ref{e2.6}) of $f$, then  $\widehat{I}_{i_{1},\dots, i_{k}}=c_{i_{1},\dots, i_{k}}\circ\widehat E$.
\begin{Th}\label{te2.17}
\begin{itemize}
\item[(1)]
An element $g\in G_{f}(X)$ belongs to $\overline{G_{f}(X)}_{n}$ if and only if all basic iterated integrals of order $\leq n$ vanish at $g$.
\item[(2)]
$$
\overline{G(X)}_{n}=\overline{G_{f}(X)}_{n}\cap G(X).
$$
\end{itemize}
\end{Th}
Let us recall that the order of a basic iterated integral is the number of its indices. \\
{\bf Proof.} We use the notation of section 2.4.2. Since the map
$q_{l}: G(X)\to Q_{l}$, see (\ref{e2.16}), is surjective, 
$q_{l}(G(X)_{n})=(Q_{l})_{n}$. Moreover, since $Q_{l}$ is a complex nilpotent Lie group, $(Q_{l})_{n}$ is a complex nilpotent Lie subgroup of $Q_{l}$. In particular, $q_{l}(\overline{G(X)}_{n})=(Q_{l})_{n}$ (because
$q_{l}(G(X)_{n})$ is dense in $q_{l}(\overline{G(X)}_{n})$). Similarly,
$q_{l}(\overline{G_{f}(X)}_{n})=q_{l}(G_{f}(X)_{n})=(Q_{l})_{n}$.
Further, the Lie algebra ${\cal L}_{Lie}^{l}$ of $Q_{l}$ is the free nilpotent Lie subalgebra of ${\cal L}_{G_{l}}$ generated by elements $\overline{X}_{1},\dots, \overline{X}_{l-1}$, see (\ref{e2.21}). Then the $n$th term $({\cal L}_{Lie}^{l})_{n}$ of the lower central series of
${\cal L}_{Lie}^{l}$ consists of the elements of form (\ref{e2.21}) with
all $c_{i_{1},\dots,i_{k}}=0$ for $k\leq n$ (i.e., the number of brackets in each term of this formula must be $\geq n$). Since the exponential map $\exp$ maps $({\cal L}_{Lie}^{l})_{n}$ surjectively onto $(Q_{l})_{n}$,
an explicit computation shows that $(Q_{l})_{n}\subset Q_{l}$ consists of elements of form (\ref{e2.15}) with $c_{0}=1$ and all $c_{i_{1},\dots, i_{k}}=0$ for $k\leq n$. 

Now, suppose that $g\in \overline{G_{f}(X)}_{n}$. Identifying $G_{f}(X)$ with its image under $\widehat E$ we have
\begin{equation}\label{e2.23}
g=I+\sum_{n=1}^{\infty}\left(\ \!\sum_{i_{1}+\cdots+ i_{k}=n}c_{i_{1},\dots, i_{k}}(g)\ \!X_{i_{1}}\cdots X_{i_{k}}\right) t^{n}.
\end{equation}
Since $q_{l}(g)\in (Q_{l})_{n}$ for any $l$, formula (\ref{e2.15}) and
the above description of $(Q_{l})_{n}$ imply that $c_{i_{1},\dots, i_{k}}(g)=0$ for all $k\leq n$. Equivalently, all basic iterated integrals of order $\leq n$ vanish at $g$. 

Conversely, assume that $g\in G_{f}(X)$ is of form (\ref{e2.23}) with $c_{i_{1},\dots, i_{k}}(g)=0$ for all $k\leq n$. Then from the description of $(Q_{l})_{n}$ and (\ref{e2.15}) we obtain that $q_{l}(g)\in (Q_{l})_{n}$ for any $l\geq 2$. Since$(Q_{l})_{n}=q_{l}(G(X)_{n})$, there are elements
$g_{l}\in G(X)_{n}$, $l\geq 2$, such that $q_{l}(g^{-1}g_{l})=I$. As in the proof of Theorem \ref{te2.10} this implies that $\{g_{l}\}_{l\geq 2}$ converges to $g$ in the topology of $G_{f}(X)$,
that is, $g\in\overline{G_{f}(X)}_{n}$. This proves (1).

Further, if $g\in \overline{G_{f}(X)}_{n}\cap G(X)$, then as above there is a  sequence $\{g\}_{l\geq 2}$ with $g_{l}\in G(X)_{n}$ such that $\lim_{l\to\infty} g_{l}=g$. Thus $g\in\overline{G(X)}_{n}$. Since the implication $\overline{G(X)}_{n}\subset\overline{G_{f}(X)}_{n}\cap G(X)$
is obvious, we obtain the proof of (2).\ \ \ \ \ $\Box$\\

{\bf 2.5.2.} By $X_{*}\subset X$ we denote the set of elements $a\in X$ such that $I_{s}(a)=0$ for all $s\in\N$, i.e., $a\in X_{*}$ if and only if its first integral $\widetilde a:I_{T}\to\Co^{\infty}$ is a closed path. We call the image $G(X_{*}):=\pi(X_{*})\subset G(X)$ the {\em subgroup of closed paths} and its closure $G_{f}(X_{*})$ in $G_{f}(X)$ the {\em subgroup of formal closed paths}. According to Theorem \ref{te2.17}, $G(X_{*})=\overline{G(X)}_{1}$ and $G_{f}(X_{*})=\overline{G_{f}(X)}_{1}$.
In this section we describe the topological lower central series of $G(X_{*})$ and $G_{f}(X_{*})$. 

Given $a=(a_{1},a_{2},\dots)\in X$ we define
\begin{equation}\label{e2.24}
\widetilde a_{i}(x):=\int_{0}^{x}a_{i}(s)\ \!ds,\ \ \ x\in I_{T}.
\end{equation}
By ${\cal P}(X)$ we denote the set of functions on $X\times I_{T}$ of the form
\begin{equation}\label{e2.25}
(\widetilde a_{i_{1}}(x))^{n_{1}}\cdots (\widetilde a_{i_{k}}(x))^{n_{k}}\cdot a_{i_{k+1}}(x),\ \ \
i_{1},\dots, i_{k+1}\in\N,\ n_{1},\dots, n_{k}\in\Z_{+}.
\end{equation}
\begin{D}\label{d2.18}
A moment of order $k$ on $X$ is an iterated integral of the form
\begin{equation}\label{e2.26}
m(a):=\int\cdots\int_{0\leq s_{1}\leq\cdots\leq s_{k}\leq T}p_{k}(a;s_{k})\cdots p_{1}(a;s_{1})\ \!ds_{k}\cdots ds_{1}
\end{equation}
where each $p_{j}\in {\cal P}(X)$, $1\leq j\leq k$.
\end{D}

Moments of the first order play an important role in the study of the center problem for Abel differential equations (see, e.g., [AL], [BFY1], [BFY2], [Y]). Also, it was proved in [Br4, Theorem 2.1] that such moments determine centers of equations (\ref{e1}) whose coefficients are either polynomials in $e^{\pm 2\pi ix/T}$ or in $x$, a result on complexity of the set of centers for these equations. (E.g., this class contains equations obtained from the Poincar\'{e} Center-Focus problem, see the Introduction.)

According to the Ree shuffle formula (\ref{e2.22}) each moment $m$ is a linear combination with natural coefficients of some basic iterated integrals. In particular, there is a continuous function $\widehat m$ on $G(X)$ from the vector space generated by all basic iterated integrals on $G(X)$ such that $m=\widehat m\circ\pi$. Thus, every such $\widehat m$ admits a continuous extension (denoted by the same symbol) to $G_{f}(X)$.  The extended function will be called a {\em moment} on $G_{f}(X)$.
By the definition the {\em order} of $\widehat m$ is the order of the moment $m$ on $X$ representing $\widehat m$.
\begin{Th}\label{te2.19}
\begin{itemize}
\item[(1)]
An element $g\in G_{f}(X_{*})$ belongs to $\overline{G_{f}(X_{*})}_{n}$ if and only if all moments of order $\leq n$ vanish at $g$.
\item[(2)]
$$
\overline{G(X_{*})}_{n}=\overline{G_{f}(X_{*})}_{n}\cap G(X_{*}).
$$
\end{itemize}
\end{Th}
{\bf Proof.} We first prove the particular case of (1) for $g\in G(X_{*})$.
Namely we will prove
\begin{Proposition}\label{pr2.20}
An element $g\in G(X_{*})$ belongs to $\overline{G(X_{*})}_{n}$ if and only if all moments of order $\leq n$ vanish at $g$.
\end{Proposition}
{\bf Proof.} This result was stated in [Br5, Theorem 3.2]. In its proof given in [Br5] some details were omitted. Here we will give
the complete proof of this fact.

First, assume that $g\in G(X_{*})_{n}$. Since $g$ represents a closed path in $\Co^{\infty}$ all moments of order $\leq n$ vanish at $g$ (see, e.g.,
[H] for properties of iterated integrals over closed paths). Since each moment is a continuous function on $G(X_{*})$, by continuity we obtain also that for $g\in\overline{G(X_{*})}_{n}$ all moments of order $\leq n$ vanish at $g$. Thus we must prove a converse statement.

So assume that $g\in G(X_{*})$ is such that all moments of order $\leq n$ vanish at $g$. We will prove that $g\in\overline{G(X_{*})}_{n}$.

For an element $p=\widetilde a_{i_{1}}^{n_{1}}\cdots\widetilde a_{i_{k}}^{n_{k}}\cdot a_{i_{k+1}}\in {\cal P}(X)$ the number
$i_{1}n_{1}+\dots + i_{k}n_{k}+i_{k+1}$ will be called the degree of $p$.
Next, for a moment $m$ on $X$ its degree $deg(m)$ is the maximum of degrees of elements $p_{j}\in {\cal P}(X)$ in its definition, see (\ref{e2.26}). In turn, the degree of the moment $\widehat m$ on $G(X)$ representing $m$ is defined as $deg(\widehat m):=deg(m)$.

We retain the notation of section 2.4.2. Also, for $Q_{l}:=q_{l}(G(X))\subset G_{l}$ we set $R_{l}:=[Q_{l},Q_{l}]$. Our proof is based on the following
\begin{Lm}\label{le2.21}
Suppose that $g\in G(X_{*})$ is such that all moments of order $\leq n$ and of degree $\leq l-1$ vanish at $g$. Then $q_{l}(g)\in (R_{l})_{n}$.
\end{Lm}
{\bf Proof.} Let $a=(a_{1},a_{2},\dots)\in X$ be such that $\pi(a)=g$.
By the definitions of $E$, see section 2.3.2, and of $q_{l}$, see (\ref{e2.15}), $q_{l}(g)$ is the monodromy of the equation
\begin{equation}\label{e2.27}
H'(x)=\left(\ \!\sum_{i=1}^{l-1}a_{i}(x)\ \!\overline{X}_{i}\right)H(x),\ \ \ x\in I_{T}.
\end{equation}
This equation can be solved by Picard iteration to obtain a solution
$H_{a}:I_{T}\to G_{l}$, $H_{a}(0)=I$, whose coefficients in expansion in $\overline{X}_{1},\dots, \overline{X}_{l-1}$ are Lipschitz functions on $I_{T}$. Then $q_{l}(g):=H_{a}(T)$. We write
\begin{equation}\label{e2.28}
H_{a}=H_{1}\cdots H_{l-1}\cdot H_{l}\ \ \ {\rm where}\ \ \ 
H_{i}:=e^{\widetilde a_{i}\overline{X}_{i}},\ \ \ 1\leq i\leq l-1,
\end{equation}
see (\ref{e2.24}). Since $g\in G(X_{*})$, $H_{i}(T)=I$ for $1\leq i\leq l-1$. This implies that
$q_{l}(g)=H_{a}(T)=H_{l}(T)$. From (\ref{e2.28}) follows that $H_{l}$ satisfies the equation
\begin{equation}\label{e2.29}
\begin{array}{c}
\displaystyle
H_{l}'=\omega\cdot H_{l}\ \ \ \ \ {\rm where}\\
\\
\displaystyle
\omega:=F^{-1}\cdot\left(\ \!\sum_{i=1}^{l-1}a_{i}\ \!\overline{X}_{i}\right)\cdot F-F^{-1}\cdot F',\ \ \ F:=H_{1}\cdots H_{l-1}.
\end{array}
\end{equation}
{\bf Claim.}\quad {\em $\omega$ is a function on $I_{T}$ with values in
the Lie algebra ${\cal L}_{R_{l}}$ of $R_{l}$.}

Indeed, the first term in the definition of $\omega$ is the logarithm of
$$
F^{-1}(x)\cdot\exp\left(\ \!\sum_{i=1}^{l-1}a_{i}(x)\overline{X}_{i}\right)\cdot F(x),\ \ \ x\in I_{T}.
$$
By the definition of $F$ for any $x\in I_{T}$ each term of this product belongs to $Q_{l}$ (see section 2.4.2 after formula (\ref{e2.21})).
Thus its logarithm belongs to the Lie algebra ${\cal L}_{Q_{l}}$ of $Q_{l}$.
Next, the second term in the definition of $\omega$ is equal to
$$
\left(\ \!\sum_{s=2}^{l-1}\ \!H_{l-1}^{-1}(x)\cdots H_{s}^{-1}(x)\cdot \widetilde a_{s-1}(x)\ \!\overline{X}_{s-1}\cdot H_{l-1}(x)\cdots H_{s}(x)\right)+\widetilde a_{l-1}(x)\ \!\overline{X}_{l-1},\ \ \ x\in I_{T}.
$$
By the same reason as above, for any $x\in I_{T}$ each term of this sum belongs to ${\cal L}_{Q_{l}}$. Thus $\omega(x)\in {\cal L}_{Q_{l}}$ for any $x\in I_{T}$. Observe also that from (\ref{e2.29}) follows that $\omega(x)$ considered as a polynomial in $\overline{X}_{i}$, $1\leq i\leq l-1$, does not contain linear terms. Then by the definition of $R_{l}$, see also the proof of Theorem \ref{te2.17}, $\omega(x)\in {\cal L}_{R_{l}}$.

This completes the proof of the Claim.

An explicit computation of $\omega$ based on the Campbell-Hausdorff formula
shows that
\begin{equation}\label{e2.30}
\omega(x)=\sum_{n=1}^{l-1}\left(\ \!\sum_{i_{1}+\dots +i_{k}=n}c_{i_{1},\dots, i_{k}}(x)\ \![\overline{X}_{i_{1}},[\overline{X}_{i_{2}},[\ \cdots , [\overline{X}_{i_{k-1}},\overline{X}_{i_{k}}]\cdots \ ]]]\right)
\end{equation}
where each $c_{i_{1},\dots, i_{k}}\in span_{\Q}({\cal P}(X))$. Moreover,
each term of ${\cal P}(X)$ in the definition of $c_{i_{1},\dots, i_{k}}$ is of degree $i_{1}+\cdots +i_{k}$.

We set 
$$
S:=\{\ \![\overline{X}_{i_{1}},[\overline{X}_{i_{2}},[\ \cdots , [\overline{X}_{i_{k-1}},\overline{X}_{i_{k}}]\cdots \ ]]]\neq 0\ :\ i_{1}+\cdots +i_{k}=n,\ 1\leq n\leq l-1\},
$$
and arrange $S$ into a sequence $\{v_{i}\}_{1\leq i\leq N}$, $N:=\# S$. Then $\omega$ can be written as
$$
\omega(x)=\sum_{i=1}^{N}f_{i}(x)\ \!v_{i},\ \ \ x\in I_{T},
$$
where $f_{i}\in span_{\Q}({\cal P}(X))$ and each term of ${\cal P}(X)$ in the definition of $f_{i}$ is of degree $\leq l-1$.

Let $\Co[[X_{1},\dots, X_{N}]]$ be the associative algebra with unit $I$ of noncommutative formal power series in free variables $X_{1},\dots, X_{N}$.
Then there is a homomorphism $\phi:\Co[[X_{1},\dots, X_{N}]]\to A$, where $A$ is the associative complex subalgebra of ${\cal A}_{l}$ (see section 2.4.2) generated by $I$ and $v_{1},\dots, v_{n}$, determined by $\phi(X_{i})=v_{i}$,
$1\leq i\leq N$, $\phi(I)=I$.

Next, consider the equation
$$
G'(x)=\left(\ \!\sum_{i=1}^{N}f_{i}(x) X_{i}\right) G(x),\ \ \ x\in I_{T}.
$$
Solving this equation by Picard iteration we get a solution
$G: I_{T}\to\Co[[X_{1},\dots, X_{N}]]$, $G(0)=I$, whose coefficients in expansion in $X_{1},\dots, X_{N}$ are Lipschitz functions on $I_{T}$.
Also, by the definition we have $\phi(G(T))=H_{l}(T)=q_{l}(g)$. Observe now
that
$$
G(T)=I+\sum_{k=1}^{\infty}\left(\ \!\sum_{1\leq i_{1},\dots, i_{k}\leq N}I_{i_{1},\dots, i_{k}}(f)\ \!X_{i_{1}}\cdots X_{i_{k}}\right),
$$
cf. (\ref{e2.11}), where $f=(f_{1},\dots, f_{N},0,\dots)\in X$. By the definition, each basic iterated integral $I_{i_{1},\dots, i_{k}}(f)$ in this formula is a linear combination of moments of order $k$ and of degree $\leq l-1$ on $X$.
In particular, by the hypothesis of the lemma, $I_{i_{1},\dots, i_{k}}(f)=0$ for all $k\leq n$, i.e., the first sum in the definition of $G(T)$ can be considered for $k\geq n+1$ only. Further, by the Ree theorem [R] $\log(G(T))$ is a Lie element in $\Co[[X_{1},\dots, X_{N}]]$. Since the degree of each monomial in $G(T)$ is greater than $n$, 
$$
\log(G(T))=\sum_{k=n+1}^{\infty}\left(\sum_{1\leq i_{1},\dots , i_{k}\leq N}g_{i_{1},\dots , i_{k}}[X_{i_{1}},[X_{i_{2}},[\ \cdots , [X_{i_{k-1}},X_{i_{k}}]\cdots \ ]]]\right).
$$
This implies that 
$$
\begin{array}{c}
\displaystyle
\phi(\log(G(T)))=\log(\phi(G(T)))=\log(q_{l}(g))=\\
\\
\displaystyle
\sum_{k=n+1}^{\infty}\left(\sum_{1\leq i_{1},\dots , i_{k}\leq N}g_{i_{1},\dots , i_{k}}[v_{i_{1}},[v_{i_{2}},[\ \cdots , [v_{i_{k-1}},v_{i_{k}}]\cdots \ ]]]\right),
\end{array}
$$
that is, $\log(q_{l}(g))\in ({\cal L}_{R_{l}})_{n}$. Since 
the exponential map $\exp$ maps $({\cal L}_{R_{l}})_{n}$ surjectively onto
$(R_{l})_{n}$, $q_{l}(g)\in (R_{l})_{n}$.

The proof of the lemma is complete.\ \ \ \ \ $\Box$

Let us finish the proof of the proposition. 

For an element $g\in G(X_{*})$ such that all moments of order $\leq n$ vanish at $g$ by Lemma \ref{le2.21} we have $q_{l}(g)\in (R_{l})_{n}$ for all $l\geq 2$. Since $q_{l}$ maps $G(X_{*})$ surjectively onto $R_{l}$ (see the proof of Theorem \ref{te2.10}), there are elements $g_{l}\in G(X_{*})_{n}$, $l\geq 2$, such that $q_{l}(g^{-1}\cdot g_{l})=I$ for all $l$. As in the proof of Theorem \ref{te2.17} this implies that $\lim_{l\to\infty} g_{l}=g$. Thus
$g\in\overline{G(X_{*})}_{n}$. This completes the proof of the proposition.
\ \ \ \ \ $\Box$

Using Proposition \ref{pr2.20} we prove now Theorem \ref{te2.19}. 

Assume that $g\in\overline{G_{f}(X_{*})}_{n}$.  Since $G_{f}(X_{*})$ is the closure in $G_{f}(X)$ of $G(X_{*})$, $\overline{G_{f}(X_{*})}_{n}$ is the closure in $G_{f}(X)$ of $G(X_{*})_{n}$. In particular,
all moments of order $\leq n$ vanish at $g$ (see the beginning of the proof of Proposition \ref{pr2.20}).

Conversely, suppose that $g\in G_{f}(X_{*})$ is such that all moments of order $\leq n$ vanish at $g$. We will prove first that $q_{l}(g)\in (R_{l})_{n}$, for all $l\geq 2$. 

By the Ree shuffle formula, given $l\geq 2$ each moment $\widehat m$ of order $\leq n$ and of degree $\leq l-1$ on $G_{f}(X_{*})$ can be presented as a linear combination with natural coefficients of basic iterated integrals $\widehat I_{i_{1},\dots, i_{k}}$ (on $G_{f}(X)$) of order $k\leq l+n-1$ with $1\leq i_{1},\dots, i_{k}\leq l-1$. We set $s:=(l+n-1)\cdot (l-1)+1$ and consider $Q_{s}:=q_{s}(G(X))=q_{s}(G_{f}(X))$ (cf. the proof of Proposition \ref{pr2.11}). By the definitions of $\widehat E$ and $q_{s}$, see (\ref{e2.12}) and (\ref{e2.15}), we have
$$
q_{s}(g)=I+\sum_{m=1}^{s-1}\left(\ \!\sum_{i_{1}+\cdots +i_{k}=m}\widehat{I}_{i_{1},\dots, i_{k}}(g)\ \!\overline{X}_{i_{1}}\cdots\overline{X}_{i_{k}}\right).
$$
Since $R_{s}:=[Q_{s},Q_{s}]=q_{s}(G(X_{*}))=q_{s}(G_{f}(X_{*}))$ and $g\in G_{f}(X_{*})$,  there is $\widetilde g\in G(X_{*})$ such that
$q_{s}(g)=q_{s}(\widetilde g)$. In particular, 
$$
\widehat{I}_{i_{1},\dots, i_{k}}(g)=\widehat{I}_{i_{1},\dots, i_{k}}(\widetilde g)\ \ \ {\rm for\ all}\ \ \ i_{1}+\dots +i_{k}=m,\  \ 1\leq m\leq s-1.
$$
From this by the above description of moments of order $\leq n$ and of degree $\leq l-1$ we obtain that for each such a moment $\widehat m$ on $G_{f}(X_{*})$, $\widehat m(g)=\widehat m(\widetilde g)$. Since $\widehat m(g)=0$ by our hypothesis, $\widehat m(\widetilde g)=0$, as well. Hence $\widetilde g$ satisfies the conditions of Lemma \ref{le2.21}. According to this lemma, $q_{l}(\widetilde g)\in (R_{l})_{n}$. But by our construction
$q_{l}(g)=q_{l}(\widetilde g)$. That is, $q_{l}(g)\in (R_{l})_{n}$. 
Using this and arguing as in the proof of Proposition \ref{pr2.20} we find
a sequence $\{g_{l}\}_{l\geq 2}\subset G(X_{*})_{n}$ such that
$\lim_{l\to\infty}g_{l}=g$. In particular, $g\in\overline{G_{f}(X_{*})}_{n}$.
This completes the proof of part (1) of Theorem \ref{te2.19}.
The second statement of this theorem follows from the first one and from Proposition \ref{pr2.20}.\ \ \ \ \ $\Box$

In conclusion let us mention that in [Br5, section 3.3] a topological characterization of paths representing elements of $\overline{G(X_{*})}_{n}$
is given, similar to that for elements of the set of universal centers ${\cal U}$, cf. section 2.2. 
%================
\subsect{\hspace*{-1em}. Subgroups of the Group of Formal Paths}

\quad {\bf 2.6.1.} By $X^{k}\subset X$ we denote the subset of elements
$a=(a_{1},\dots, a_{k},0,0,\dots)\in X$. Then $X^{k}$ is a sub-semigroup of $X$. The first integrals of elements of $X^{k}$ are paths in $\Co^{k}$.
We set $G(X^{k}):=\pi(X^{k})\subset G(X)$ and let $G_{f}(X^{k})$ be the closure of $G(X^{k})$ in $G_{f}(X)$. The group $G_{f}(X^{k})$ will be called the {\em group of formal paths} in $\Co^{k}$. Let $p_{k}:X\to X^{k}$,
$p_{k}(a_{1},a_{2},\dots):=(a_{1},\dots, a_{k}, 0,0,\dots)$, be the natural projection. It induces a surjective homomorphism of topological groups
$\widehat p_{k}:G_{f}(X)\to G_{f}(X^{k})$. In particular, $G_{f}(X)$ is the semidirect product of groups $Ker\ \!\widehat p_{k}$ and $G_{f}(X^{k})$.
In turn, the homomorphism $\widehat p_{k}$ determines a continuous Lie algebra homomorphism $\phi_{k}:{\cal L}_{Lie}\to {\cal L}_{Lie}$ where ${\cal L}_{Lie}$ is the Lie algebra of $G_{f}(X)$, see section 2.4.1. It is determined by the conditions
$$
\phi_{k}(X_{s}):=
\left\{
\begin{array}{ccc}
X_{s}&{\rm if}&1\leq s\leq k\\
\\
0&{\rm if}&s>k.
\end{array}
\right.
$$
The image of $\phi_{k}$ is a closed Lie subalgebra ${\cal L}_{Lie}^{k}$ of ${\cal L}_{Lie}$ consisting of Lie elements in variables $X_{1},\dots, X_{k}$ and $t$. Identifying $G_{f}(X)$ with its image under map $\widehat E$, see section 2.3.2, we obtain the commutative diagram
\begin{equation}\label{e2.31}
\begin{array}{ccc}
\ \ \ {\cal L}_{Lie}&\stackrel{\exp}{\longrightarrow}&\!\!\!G_{f}(X)\\
\phi_{k}\downarrow& &\ \downarrow \widehat p_{k}\\
\ \ \ \ \ {\cal L}_{Lie}^{k}&\stackrel{\exp}{\longrightarrow}& G_{f}(X^{k}).
\end{array}
\end{equation}
Thus ${\cal L}_{Lie}^{k}$ can be regarded as the Lie algebra of $G_{f}(X^{k})$. 

Also, analogs of Theorems \ref{te2.17} and \ref{te2.19} are valid for  topological lower central series of $G_{f}(X^{k})$ and $G_{f}(X_{*}^{k})$ (the {\em group of formal closed paths in} $\Co^{k}$) where in these results we consider basic iterated integrals and moments on $G_{f}(X^{k})$, respectively. \\

{\bf 2.6.2.} Let $\F\subset\Co$ be a field. By $X_{\F}\subset X$ we denote the subset of elements $a\in X$ such that $I(a)\in\F$ for all basic iterated integrals on $X$. Formulas (\ref{e2.2}) and (\ref{e2.3}) imply that $X_{\F}$ is a sub-semigroup of $X$. By $G(X_{\F}):=\pi(X_{\F})$ we denote the subgroup of $G_{f}(X)$ generated by $X_{\F}$. The homomorphism $\widehat E$,
see (\ref{e2.12}), embeds $G(X_{\F})$ into the subalgebra ${\cal A}_{\F}$ of ${\cal A}$, see section 2.3.1, of formal power series whose coefficients in expansion in $I, X_{1},X_{2},\dots$ and $t$ belong to $\F$. We will identify $G(X_{\F})$ with its image under $\widehat E$.

Next, by $J_{\F}\subset {\cal A}_{\F}$ we denote the two-sided ideal of elements $f$
whose series expansions do not contain terms with $I$. By $J_{\F}^{k}$ we denote the $k\ \!$th power of $J_{\F}$. Let us introduce 
the $J_{\F}$-adic topology on ${\cal A}_{\F}$, i.e., a sequence $\{f_{i}\}_{i\in\N}\subset {\cal A}_{\F}$ converges in this topology to $f\in {\cal A}_{\F}$ if and only if for any $l\in\N$ there is a natural number $N_{l}$ such that for all $n\geq N_{l}$ the images of $f_{n}$ and $f$ in the quotient algebra ${\cal A}_{\F}/J_{\F}^{l}$ coincide. Observe that ${\cal A}_{\F}$ is complete in this topology. By $G_{f}(X_{\F})\subset {\cal A}_{\F}$
we denote the completion of $G(X_{\F})$ in the $J_{\F}$-adic topology. We call it the {\em group of formal paths over} $\F$.

Let $[X_{\F}]_{rect}$ be a sub-semigroup of the semigroup of rectangular paths $X_{rect}$ generated by elements $a_{i}=(a_{1i},a_{2i},\dots)$ where
$a_{ki}=0$ for $k\neq i$ and $a_{ii}=c_{i}/T$, $c_{i}\in\F$. Then $G([X_{\F}]_{rect})$ is the subgroup of $G(X_{rect})$ generated by elements
$e^{c_{i}X_{i}t^{i}}$, $c_{i}\in\F$, $i\in\N$. Based on the results of sections 2.4.1 and 2.4.2 one obtains that $G([X_{\F}]_{rect})$ is dense in
$G_{f}(X_{\F})$.  Moreover, the Lie algebra ${\cal L}_{Lie(\F)}\subset {\cal L}_{Lie}$ of $G_{f}(X_{\F})$ consists of all Lie elements with coefficients
from $\F$. 

In the same way one can formulate analogs of Theorems \ref{te2.17} and \ref{te2.19} for topological lower central series of $G(X_{\F})$ and $G_{f}(X_{\F})$ in terms of basic iterated integrals and moments on $G_{f}(X_{\F})$, respectively. (We leave the details to the reader.)\\

{\bf 2.6.3.} In the sequel we will use the following result.

Let $R\subset {\cal L}_{Lie}$ be a subset.  By $A_{R}\subset {\cal L}_{Lie}$ we denote the minimal closed Lie subalgebra containing $R$. Consider the subgroup $H_{R}\subset G_{f}(X)\ \!(\subset G)$ generated by elements $e^{cr}$ for all possible $r\in R$ and $c\in\Co$. By $\overline{H}_{R}$ we denote the closure of $H_{R}$ in $G_{f}(X)$. 
\begin{Proposition}\label{pr2.22}
We have 
$$
\log(\overline{H}_{R})= A_{R}.
$$
Moreover, $\overline{H}_{R}$ is a normal subgroup of $G_{f}(X)$ if and only if
$A_{R}$ is a normal Lie subalgebra of ${\cal L}_{Lie}$.
\end{Proposition}
{\bf Proof.} We retain the notations of section 2.4.2. 

Consider the image
$(H_{R})_{l}:=q_{l}(H_{R})\subset Q_{l}:=q_{l}(G_{f}(X))$. Then, 
$(H_{R})_{l}$ is a complex Lie subgroup of $Q_{l}$. It can be shown similarly to the statement of Lemma \ref{le2.12} where instead of the map $F_{K}$ given by (\ref{eq2.18}) we determine now a new map $F_{K}:\Co^{k}\to\Co^{N}$, $N:=dim_{\Co}Q_{l}$, with $K=\{r_{1},\dots, r_{k}\}\subset R$ by the formula
$$
F_{K}(z_{1},\dots, z_{k}):=q_{l}(e^{z_{1}r_{1}}\cdots e^{z_{k}r_{k}}).
$$
Then, as in the proof of the lemma, $F_{K}$ is a holomorphic polynomial map.
Applying now the arguments of Lemma \ref{le2.12} to the family of such maps $F_{K}$, we finally get the required: $(H_{R})_{l}$ is a complex Lie subgroup of $Q_{l}$. In particular, we also have 
$(\overline{H}_{R})_{l}:=q_{l}(\overline{H}_{R})=(H_{R})_{l}$.

From the above statement we obtain that $\log((H_{R})_{l})\subset {\cal L}_{Lie}^{l}$ is the Lie algebra of $(H_{R})_{l}$ (cf. (\ref{e2.21})). Since $\log((H_{R})_{l})$ contains $q_{l}(R)$ and $(H_{R})_{l}$ is generated by  elements $e^{cq_{l}(r)}$ for all possible $r\in R$ and $c\in\Co$, the Campbell-Hausdorff formula implies that 
\begin{equation}\label{e2.32}
\log((\overline{H}_{R})_{l})=\log((H_{R})_{l})=q_{l}(A_{R}).
\end{equation}

Assume now that $h\in\log(\overline{H}_{R})$. According to (\ref{e2.32}) there is a sequence $\{h_{l}\}_{l\geq 2}\subset A_{R}$ such that
$q_{l}(h-h_{l})=0$ for all $l$. This implies that $\lim_{l\to\infty}h_{l}=h$, that is, $\log(\overline{H}_{R})\subset A_{R}$. The inclusion $A_{R}\subset\log(\overline{H}_{R})$ is obtained similarly using the fact that $\log(\overline{H}_{R})$ is a closed subset of ${\cal L}_{Lie}$.

Now, if $\overline{H}_{R}$ is a normal subgroup of $G_{f}(X)$, then
$(\overline{H}_{R})_{l}$ is a normal Lie subgroup of $Q_{l}$. This and (\ref{e2.32}) imply that $q_{l}(A_{R})$ is a normal Lie subalgebra of ${\cal L}_{Lie}^{l}$ (the standard fact of the theory of finite-dimensional complex Lie groups). Thus, if $a\in {\cal L}_{Lie}$, $h\in A_{R}$, then
$q_{l}([a,h)]=[q_{l}(a),q_{l}(h)]\in q_{l}(A_{R})$ for all $l$. As above, the latter implies that $\lim_{l\to\infty}g_{l}=[a,h]$ for some $\{g_{l}\}_{l\geq 2}\subset A_{R}$, i.e., $[a,h]\in A_{R}$. Hence, $A_{R}$ is a normal subalgebra of ${\cal L}_{Lie}^{l}$. The converse statement can be obtained in the same way (we leave the details to the reader).\ \ \ \ \ $\Box$
\begin{R}\label{r2.23}
{\rm The above result can be proved by means of the Campbell-Hausdorff formula only. This method of the proof works also in a similar result for $R\subset {\cal L}_{Lie(\F)}$ and $H_{R}\subset G_{f}(X_{\F})$ generated by elements $e^{cr}$ for all possible $r\in R$ and $c\in\F$.}
\end{R}
%=================
\sect{Center Problem for ODEs}

\subsect{\hspace*{-1em}. An Explicit Expression for the First Return Map}

\quad {\bf 3.1.1.} Let $\Co[[z]]$ be the algebra of formal complex power series in $z$. By $D,\ \!L:\Co[[z]]\to\Co[[z]]$ we denote the differentiation and the left translation operators defined on $f(z)=\sum_{k=0}^{\infty}c_{k}z^{k}$ by
\begin{equation}\label{e3.1}
(Df)(z):=\sum_{k=0}^{\infty}(k+1)c_{k+1}z^{k},\ \ \ (Lf)(z):=\sum_{k=0}^{\infty}c_{k+1}z^{k}.
\end{equation}
Let ${\cal A}(D,L)$ be the associative algebra with unit $I$ of complex polynomials in $I$, $D$ and $L$. By ${\cal A}(D,L)[[t]]$ we denote the associative algebra of formal power series in $t$ with coefficients from ${\cal A}(D,L)$. Also, by $G_{0}(D,L)[[t]]$ we denote the group of invertible elements of ${\cal A}(D,L)[[t]]$ consisting of elements whose expansions in $t$ begin with $I$.

Further, consider equation (\ref{e1}) corresponding to
an $a=(a_{1}, a_{2},\dots)\in X$:
\begin{equation}\label{e3.2}
\frac{dv}{dx}=\sum_{i=1}^{\infty}a_{i}(x) v^{i+1},\ \ \ x\in I_{T}.
\end{equation}
Using the linearization of (\ref{e3.2}) described in [Br1] we associate to this equation the following system of ODEs:
\begin{equation}\label{e3.3}
H'(x)=\left(\ \!\sum_{i=1}^{\infty} a_{i}(x) DL^{i-1}t^{i}\right) H(x),\ \ \ x\in I_{T}.
\end{equation}
Solving (\ref{e3.3}) by Picard iteration we obtain a solution $H_{a}: I_{T}\to G_{0}(D,L)[[t]]$, $H_{a}(0)=I$,  whose coefficients in the series expansion in $D, L$ and $t$ are Lipschitz functions on $I_{T}$. It was established in [Br2, Theorem 1.1] that (\ref{e3.2}) determines a center (i.e., $a\in {\cal C})$ if and only if $H_{a}(T)=I$. This implies the following result (see [Br5, Proposition 2.1]).
\begin{Th}\label{te3.1}
\begin{equation}\label{e3.4}
a\in {\cal C}\ \Longleftrightarrow \sum_{i_{1}+\cdots +i_{k}=i}p_{i_{1},\dots, i_{k}}I_{i_{1},\dots, i_{k}}(a)\equiv 0\ \ \ {\rm for\ all}\ \ \ i\in\N
\end{equation}
where $p_{i_{1},\dots, i_{k}}$ is a polynomial of degree $k$ defined by
\begin{equation}\label{e3.5}
p_{i_{1},\dots, i_{k}}(t)=(t-i_{1}+1)(t-i_{1}-i_{2}+1)(t-i_{1}-i_{2}-i_{3}+1)\cdots (t-i+1),\ \ \ t\in\Co.
\end{equation}
\end{Th}

{\bf 3.1.2.}
Let $G[[r]]$ be the set of formal complex power series $f(r)=r+\sum_{i=1}^{\infty}d_{i}r^{i+1}$.
Let $d_{i}:G[[r]]\to\Co$ be such that $d_{i}(f)$ is the $(i+1)$st coefficient in the series expansion of $f$. We equip $G[[r]]$ with the weakest topology in which all $d_{i}$ are continuous functions and consider the multiplication $\circ$ on $G[[r]]$ defined by the composition of series. Then $G[[r]]$ is a separable topological group. Moreover, it is contractible and residually torsion free nilpotent. By $G_{c}[[r]]\subset G[[r]]$ we denote the subgroup of power series locally convergent near $0$ equipped with the induced topology. Next, we define the map $P:X\to G[[r]]$ by the formula
\begin{equation}\label{e3.6}
P(a):=r+\sum_{i=1}^{\infty}\left(\sum_{i_{1}+\cdots +i_{k}=i}p_{i_{1},\dots, i_{k}}(i)\cdot I_{i_{1},\dots, i_{k}}(a)\right)r^{i+1},
\end{equation}
see (\ref{e3.5}). It follows from [Br3] that  $P(a*b)=P(a)\circ P(b)$ and $P(X)=G_{c}[[r]]$. 
Moreover, let $v(x;r;a)$, $x\in I_{T}$, be the Lipschitz solution of equation (\ref{e3.2}) with initial value $v(0;r;a)=r$. Clearly for every $x\in I_{T}$ we have $v(x;r;a)\in G_{c}[[r]]$. It was proved in [Br1] that $P(a)=v(T;\cdot;a)$ (i.e., $P(a)$ is the {\em first return map} of (\ref{e3.2})). In particular, we have
\begin{equation}\label{e3.7}
a\in {\cal C}\ \Longleftrightarrow\ \sum_{i_{1}+\cdots +i_{k}=i}p_{i_{1},\dots, i_{k}}(i)\cdot I_{i_{1},\dots, i_{k}}(a)\equiv 0\ \ \ {\rm for\ all}\ \ \ i\in\N.
\end{equation}

Also, (\ref{e3.6}) implies that there is a continuous homomorphism $\widehat P: G(X)\to G[[r]]$ such that $P=\widehat P\circ\pi$
(where $\pi: X\to G(X)$ is the quotient map). We extend it by continuity to $G_{f}(X)$ retaining the same symbol for the extension. 
%=========
\subsect{\hspace*{-1em}. Group of Formal Centers}
\quad {\bf 3.2.1.} Let $\Co\!\left\langle X_{1},X_{2}\right\rangle$ be the associative algebra with unit $I$ of complex polynomials in $I$ and free noncommutative variables $X_{1}$, $X_{2}$. Consider a homomorphism
$\phi:\Co\!\left\langle X_{1},X_{2}\right\rangle\to {\cal A}(D,L)$  defined by conditions: $\phi(X_{1}):=D$, $\phi(X_{2}):=L$. 
Then $Ker\phi\subset \Co\!\left\langle X_{1},X_{2}\right\rangle$ is a two-sided ideal generated by the element
$X_{1}X_{2}-X_{2}X_{1}+X_{2}^{2}$, see [Br2, Proposition 2.2]. In particular, 
see [Br2, Lemma 2.4], each $p\in {\cal A}(D,L)$ is uniquely presented as
\begin{equation}\label{e3.8}
p(D,L,I)=a_{0}I+\sum_{k=1}^{n}F_{k}(D,L)\ \ \ {\rm where}\ \ \
F_{k}(D,L)=\sum_{i=0}^{k}a_{i k-i, k}D^{i}L^{k-i}
\end{equation}
with all $a_{0}, a_{ij,k}\in\Co$. 

We say that such $p$ has {\em degree} $n$ if the polynomial $p(x,y,1)$ in commutative variables $x,y$ has degree $n$. By $P_{n}\subset {\cal A}(D,L)$ we denote the complex vector space of polynomials of degree $\leq n$. We naturally identify $P_{n}$ with the hermitian space $\Co^{k(n)}$ where $k(n)=dim_{\Co}P_{n}$.

Let ${\cal A}_{*}\subset {\cal A}(D,L)[[t]]$ be the subalgebra of series
\begin{equation}\label{e3.9}
f=\sum_{n=0}^{\infty}f_{n}t^{n}\ \ \ {\rm with}\ \ \ f_{n}\in P_{n},\ \ \ n\in\Z_{+}.
\end{equation}
We equip ${\cal A}_{*}$ with the weakest topology in which all coefficients $f_{n}$ in expansion (\ref{e3.9}) considered as functions in $f$ are continuous maps
of ${\cal A}_{*}$ into $\Co^{k(n)}$, $n\in\Z_{+}$. Since the set of such maps is countable, ${\cal A}_{*}$ is metrizable, cf. (\ref{e2.4}). Moreover, if $d$ is a metric on ${\cal A}_{*}$ compatible with topology, then $({\cal A}_{*},d)$ is a complete metric space (i.e., a sequence $\{f_{k}=\sum_{n=0}^{\infty}f_{nk}t^{n}\}_{k\in\N}\subset {\cal A}_{*}$ such that each $\{f_{nk}\}_{k\in\N}\subset P_{n}$ is a Cauchy sequence converges to $f=\sum_{n=0}^{\infty}f_{n}t^{n}\in {\cal A}_{*}$ where $f_{n}=\lim_{k\to\infty}f_{nk}$, $n\in\Z_{+}$).

By $G_{*}\subset G_{0}(D,L)[[t]]$ we
denote the subgroup of elements $f\in {\cal A}_{*}$ with $f_{0}=I$ equipped with the induced topology. Then $(G_{*}, d)$ is a complete metric space.

Next, consider an algebra homomorphism $\Psi :{\cal A}\to {\cal A}_{*}$ determined by conditions 
\begin{equation}\label{e3.10}
\Psi(X_{i}):=DL^{i-1},\ \ \  i\in\N.
\end{equation}
(Recall that 
${\cal A}\subset \Co\left\langle X_{1},X_{2},\dots\right\rangle [[t]]$ is defined by (\ref{e2.6}).) By the definition for
$$
f=c_{0}I+\sum_{n=1}^{\infty}\left(\ \!\sum_{i_{1}+\cdots +i_{k}=n}c_{i_{1},\dots, i_{k}}X_{i_{1}}\cdots X_{i_{k}}\right)t^{n}\in {\cal A}
$$
we have
$$
\Psi(f):=c_{0}I+\sum_{n=1}^{\infty}\left(\ \!\sum_{i_{1}+\cdots +i_{k}=n}c_{i_{1},\dots, i_{k}}DL^{i_{1}-1}\cdots DL^{i_{k}-1}\right)t^{n}\in {\cal A}_{*}.
$$
Rewriting each $DL^{i_{1}-1}\cdots DL^{i_{k}-1}$ in the form (\ref{e3.8}) 
(using identity $DL-LD=-L^{2}$), we conclude that $\Psi$ is a continuous homomorphism of topological algebras.
Moreover,
$\Psi|_{G}:G\to G_{*}$ is a continuous homomorphism of topological groups. (Recall that $G$ is the subset of elements of ${\cal A}$ whose expansions in $t$ begin with $I$.) 

Observe that $\Psi$ transfers equation (\ref{e2.8}) (determining  $E: X\to G$, see (\ref{e2.9})) to equation (\ref{e3.3}). In particular,
we have 
$$
\Psi(E(a))=H_{a},\ \ \ a\in X.
$$
The last identity gives rise to the formula
\begin{equation}\label{e3.11}
\Psi(\widehat E(g))=I+\sum_{n=1}^{\infty}\left(\ \!\sum_{i_{1}+\cdots +i_{k}=n}\widehat I_{i_{1},\dots, i_{k}}(g)DL^{i_{1}-1}\cdots DL^{i_{k}-1}\right) t^{n},\ \ \ g\in G_{f}(X).
\end{equation}
(Here, as before, we regard the basic iterated integrals $\widehat I_{\cdot}$ as continuous functions on $G_{f}(X)$ extending them by continuity from $G(X)$.)\\

{\bf 3.2.2.} Let us observe that the Lie algebra ${\cal L}_{G_{*}}$ of $G_{*}$ consists of elements of ${\cal A}_{*}$ of form (\ref{e3.9}) with $f_{0}=0$. As usual, for $f,g\in {\cal L}_{G_{*}}$ their product is defined by the formula $[f,g]:=f\cdot g-g\cdot f$. Also, the map $\exp: {\cal L}_{G_{*}}\to G_{*}$, $\exp(f):=e^{f}$, is a homeomorphism.

The group homomorphism $\Psi: G\to G_{*}$ determines a continuous homomorphism of the corresponding Lie algebras such that the following diagram is commutative:
\begin{equation}\label{e3.12}
\begin{array}{ccc}
\ \ \ {\cal L}_{G}&\stackrel{\exp}{\longrightarrow}&\!\!\!G\\
\Psi\downarrow& &\ \downarrow \Psi\\
\ \ \ \ \ {\cal L}_{G_{*}}&\stackrel{\exp}{\longrightarrow}& G_{*}.
\end{array}
\end{equation}
By ${\cal L}_{S}\subset {\cal L}_{G_{*}}$ we denote the image under $\Psi$ of the Lie algebra ${\cal L}_{Lie}$ of the group $\widehat E(G_{f}(X))\ \! (\cong G_{f}(X))$, see section 2.4.1. 

In the calculations below we will use the following result.
\begin{Lm}\label{le3.2}
$$
[DL^{i},DL^{j}]=(i-j)DL^{i+j+1}.
$$
\end{Lm}
{\bf Proof.} It suffices to check the identity for elements  $z^{n}\in\Co[[z]]$ with $n\geq i+j+2$. Then we have
$$
\begin{array}{c}
\displaystyle
(DL^{i}DL^{j})(z^{n})=(n-j-i-1)(n-j)z^{n-j-i-2},\\
\\
\displaystyle
(DL^{j}DL^{i})(z^{n})=(n-i-j-1)(n-i)z^{n-i-j-2},\ \ \ {\rm and}\\
\\
\displaystyle
(i-j)(DL^{i+j+1})(z^{n})=(i-j)(n-i-j-1)z^{n-i-j-2}.
\end{array}
$$
These identities imply the required result.\ \ \ \ \ $\Box$

Now, for an element
$$
g=\sum_{n=1}^{\infty}\left(\ \!\sum_{i_{1}+\dots +i_{k}=n}c_{i_{1},\dots, i_{k}}[X_{i_{1}},[X_{i_{2}},[\ \cdots , [X_{i_{k-1}},X_{i_{k}}]\cdots \ ]]]\right) t^{n}\in {\cal L}_{Lie}
$$
we have
$$
\Psi(g):=
\sum_{n=1}^{\infty}\left(\ \!\sum_{i_{1}+\dots +i_{k}=n}c_{i_{1},\dots, i_{k}}[DL^{i_{1}-1},[DL^{i_{2}-1},[\ \cdots , [DL^{i_{k-1}-1},DL^{i_{k}-1}]\cdots \ ]]]\right) t^{n}.
$$
Simplifying the right-hand side by Lemma \ref{le3.2} we finally obtain
\begin{equation}\label{e3.13}
\Psi(g)=\sum_{n=1}^{\infty}\left(\ \!\sum_{i_{1}+\cdots + i_{k}=n}c_{i_{1},\dots, i_{k}}\cdot\gamma_{i_{1},\dots, i_{k}}DL^{n-1}\right)t^{n}
\end{equation}
where $\gamma_{n}=1$ and
$$
\gamma_{i_{1},\dots, i_{k}}:=
(i_{k-1}-i_{k})(i_{k-1}+i_{k}-i_{k-2})\cdots
(i_{2}+\cdots +i_{k}-i_{1})\ \ \ {\rm for}\ \ \ k\geq 2.
$$
\begin{Proposition}\label{pr3.3}
The following is true:
\begin{equation}\label{e3.14}
{\cal L}_{S}=\left\{g\in {\cal L}_{G_{*}}\ :\ g=\sum_{n=1}^{\infty}g_{n}\ \!DL^{n-1}t^{n},\ g_{n}\in\Co,\ n\in\N\right\}.
\end{equation}
\end{Proposition}
{\bf Proof.} By $V$ we denote the vector space on the right-hand side of (\ref{e3.14}). According to Lemma \ref{le3.2}, $V$ is a closed Lie subalgebra of ${\cal L}_{G_{*}}$. Moreover, by (\ref{e3.13}) ${\cal L}_{S}\subset V$.

The converse implication is obvious:
for an element $g=\sum_{n=1}^{\infty}g_{n}DL^{n-1}t^{n}\in V$ consider $\widetilde g:=\sum_{n=1}^{\infty}g_{n}X_{n}t^{n}\in {\cal L}_{Lie}$.
Then $\Psi(\widetilde g)=g$, i.e., $g\in {\cal L}_{S}$.\ \ \ \ \ $\Box$

From Proposition \ref{pr3.3} and diagram (\ref{e3.12}) we immediately obtain:
\begin{itemize}
\item[]
{\em 
$S:=(\Psi\circ\widehat E)(G_{f}(X))$ is a closed subgroup of $G_{*}$ with the Lie algebra ${\cal L}_{S}$.}
\end{itemize}

{\bf 3.2.3.} 
According to (\ref{e3.11}) the normal subgroup $Ker(\Psi\circ\widehat E)\subset G_{f}(X)$ consists of elements $g$ such that
$$
\sum_{i_{1}+\cdots +i_{k}=n}\widehat I_{i_{1},\dots, i_{k}}(g)DL^{i_{1}-1}\cdots DL^{i_{k}-1}=0\ \ \ {\rm for\ all}\ \ \ n\in\N .
$$
Repeating literally the arguments of the proof of [Br5, Proposition 2.1] we obtain 
\begin{equation}\label{e3.15}
g\in Ker(\Psi\circ\widehat E)\ \Longleftrightarrow \sum_{i_{1}+\cdots +i_{k}=i}p_{i_{1},\dots, i_{k}}\widehat I_{i_{1},\dots, i_{k}}(g)\equiv 0\ \ \ {\rm for\ all}\ \ \ i\in\N
\end{equation}
where $p_{i_{1},\dots, i_{k}}$ is the polynomial defined in (\ref{e3.5}).

Next, recall that the homomorphism $\widehat P: G_{f}(X)\to G[[r]]$ is determined by 
\begin{equation}\label{e3.16}
\widehat P(g):=r+\sum_{i=1}^{\infty}\left(\sum_{i_{1}+\cdots +i_{k}=i}p_{i_{1},\dots, i_{k}}(i)\cdot \widehat I_{i_{1},\dots, i_{k}}(g)\right)r^{i+1},
\end{equation}
see section 3.1.2. Then from (\ref{e3.15}) we obtain that
$$
Ker(\Psi\circ\widehat E)\subset Ker\widehat P.
$$
This implies that there is a homomorphism $\Phi:S\to G[[r]]$ such that
\begin{equation}\label{e3.17}
\widehat P=\Phi\circ\Psi\circ\widehat E.
\end{equation}
\begin{Proposition}\label{pr3.4}
$\Phi: S\to G[[r]]$ is an isomorphism of topological groups.
\end{Proposition}
{\bf Proof.} Suppose that
\begin{equation}\label{e3.18}
s:=\exp\left(\ \!\sum_{n=1}^{\infty}s_{n}DL^{n-1}t^{n}\right)\in S.
\end{equation}
Then, from (\ref{e3.16}), (\ref{e3.17}) we obtain
\begin{equation}\label{e3.19}
\Phi(s)=r+\sum_{i=1}^{\infty}\left(\ \!\sum_{i_{1}+\cdots+i_{k}=i}
\frac{p_{i_{1},\dots, i_{k}}(i)\cdot s_{i_{1}}\cdots s_{i_{k}}\cdot T^{k}}{k!}\right) r^{i+1}.
\end{equation}
Since the map $\log:S\to {\cal L}_{S}$, $\log(f):=-\sum_{i=1}^{\infty}\frac{(I-f)^{i}}{i}$, is a homeomorphism, formula (\ref{e3.19}) and the definitions of topologies on $S$ and $G[[r]]$ imply that $\Phi$ is a continuous homomorphism.

Further, the expression in the brackets of (\ref{e3.19}) can be written in
the form
$s_{i}T+p_{i}(s_{1}T,\dots, s_{i-1}T)$
where $p_{i}$ is a polynomial of degree $i$ with rational coefficients on
$\Re^{i-1}$. In particular, for any sequence $\{d_{i}\}_{i\in\N}\subset\Co$
one can solve consequently the equations
$$
s_{i}T+p_{i}(s_{1}T,\dots, s_{i-1}T)=d_{i},\ \ \ i\in\N,
$$
to get an element $s$ of form (\ref{e3.18}) such that $\Phi(s)=r+\sum_{i=1}^{\infty}d_{i}r^{i+1}$. Moreover, each $s_{i}$ in the definition of $s$ is a polynomial in variables $d_{1},\dots, d_{i}$. This implies that $\Phi$ has a continuous inverse $\Phi^{-1}:G[[r]]\to S$ and completes the proof of the proposition.\ \ \ \ \ $\Box$
\begin{R}\label{re3.5}
{\rm The Lie algebra ${\cal L}_{S}$ is isomorphic to the algebra $W_{1}(1)$, the nilpotent part of the Witt algebra of formal vector fields on $\Re$, which is known to be the Lie algebra of $G[[r]]$. Recall that $W_{1}(1)$ has the natural basis $e_{i}:=r^{i+1}\frac{d}{dr}$, $i\in\N$. Then the isomorphism $w: {\cal L}_{S}\to W_{1}(1)$ is given by $w(DL^{i-1})=-e_{i}$, $i\in\N$. Identifying ${\cal L}_{S}$ with $W_{1}(1)$ by $w$ we can regard the map $\Phi\circ\exp$ as an exponential map $W_{1}(1)\to G[[r]]$.}
\end{R}

According to (\ref{e3.17}) and Proposition \ref{pr3.4} we have
$$
Ker(\Psi\circ\widehat E)=Ker\widehat P.
$$
This group is denoted by $\widehat{\cal C}_{f}$ and called the {\em group of formal centers} of equation (\ref{e1}). By the definition $\widehat{\cal C}_{f}$ is a closed normal subgroup of $G_{f}(X)$. Moreover, 
$\widehat{\cal C}_{f}$ contains the subgroup $\widehat{\cal C}:=\pi({\cal C})\subset G(X)$, the {\em group of centers} of equation (\ref{e1}).
%===============
\subsect{\hspace*{-1em}. Properties of the Group of Formal Centers}

\quad {\bf 3.3.1.} In what follows we identify $G_{f}(X)$ and $G(X)$ with their images under $\widehat E$.

Formulas (\ref{e3.15}), (\ref{e3.16}) imply, cf. section 3.1:
\begin{equation}\label{e3.20}
\begin{array}{c}
\displaystyle
g\in \widehat{\cal C}_{f}\ \ \ \Longleftrightarrow\ \sum_{i_{1}+\cdots +i_{k}=i}p_{i_{1},\dots, i_{k}}\widehat I_{i_{1},\dots, i_{k}}(g)\equiv 0\ \ \ {\rm for\ all}\ \ \ i\in\N\ \ \ \Longleftrightarrow\\
\\
\displaystyle
\sum_{i_{1}+\cdots +i_{k}=i}p_{i_{1},\dots, i_{k}}(i)\cdot \widehat I_{i_{1},\dots, i_{k}}(g)=0\ \ \ {\rm for\ all}\ \ \ i\in\N\\
\\
{\rm where}\ \ \ \ \
p_{i_{1},\dots, i_{k}}(t)=(t-i_{1}+1)(t-i_{1}-i_{2}+1)\cdots (t-i+1).
\end{array}
\end{equation}

In turn, (\ref{e3.12}) and (\ref{e3.13}) imply that the Lie algebra ${\cal L}_{\widehat{\cal C}_{f}}\subset {\cal L}_{Lie}$ of $\widehat{\cal C}_{f}$
consists of elements
$$
\sum_{n=1}^{\infty}\left(\ \!\sum_{i_{1}+\dots +i_{k}=n}c_{i_{1},\dots, i_{k}}[X_{i_{1}},[X_{i_{2}},[\ \cdots , [X_{i_{k-1}},X_{i_{k}}]\cdots \ ]]]\right) t^{n}
$$
such that
\begin{equation}\label{e3.21}
\begin{array}{c}
\displaystyle 
\sum_{i_{1}+\cdots + i_{k}=n}c_{i_{1},\dots, i_{k}}\cdot\gamma_{i_{1},\dots, i_{k}}=0\ \ \ {\rm for\ all}\ \ \ n\in\N \ \ \ {\rm where}\ \ \ \gamma_{n}=1\ \ \ {\rm and}\\
\\
\gamma_{i_{1},\dots, i_{k}}=
(i_{k-1}-i_{k})(i_{k-1}+i_{k}-i_{k-2})\cdots
(i_{2}+\cdots +i_{k}-i_{1})\ \ \ {\rm for}\ \ \ k\geq 2.
\end{array}
\end{equation}
(In particular, the map $\exp: {\cal L}_{\widehat{\cal C}_{f}}\to  \widehat{\cal C}_{f}$ is a homeomorphism.)

Further, by the definition we have 
\begin{equation}\label{e3.22}
\widehat{\cal C}=\widehat{\cal C}_{f}\cap G(X).
\end{equation}
\begin{Proposition}\label{pr3.6}
\ $\widehat{\cal C}$ is a dense subgroup of $\widehat{\cal C}_{f}$.
\end{Proposition}
{\bf Proof.} According to [Br3, Theorem 2.10] there exists a continuous embedding
$T: G_{c}[[r]]\to G(X)$ such that $\widehat P\circ T=id$. Moreover,
$\widetilde T: G_{c}[[r]]\times\widehat{\cal C}\to G(X)$, $\widetilde T(s,g):=T(s)\cdot g$ is a homeomorphism. We extend $T$ by continuity to
a map $T_{f}:G[[r]]\to G_{f}(X)$.
Since $\widehat P\circ T=id$, similarly we have $\widehat P\circ T_{f}=id$. In particular, $T_{f}$ is an embedding. 

Let $cl(\widehat{\cal C})$ be the closure of $\widehat{\cal C}$ in $\widehat{\cal C}_{f}$. We extend the map $\widetilde T$ by continuity to a map $\widetilde T_{f}:G[[r]]\times cl(\widehat{\cal C})\to G_{f}(X)$. Then, $\widetilde T_{f}(s,g):=T_{f}(s)\cdot g$. Next, since $\widetilde T$ is a homeomorphism and $G(X)$ is dense in $G_{f}(X)$, the map $\widetilde T_{f}$ is a homeomorphism, as well (the inverse $\widetilde T_{f}^{-1}$ of $\widetilde T_{f}$ is the extension by continuity of $\widetilde T^{-1}$).
In particular, for $g\in \widehat{\cal C}_{f}$ we have $g=T_{f}(s)\cdot h$ for some $s\in G[[r]]$, $h\in cl(\widehat{\cal C})$. This implies 
$$
1=\widehat P(g)=\widehat P(T_{f}(s)\cdot h)=\widehat P(T_{f}(s))\cdot\widetilde P(h)=s.
$$
Hence, $g=h\in cl(\widehat{\cal C})$.\ \ \ \ \ $\Box$\\

{\bf 3.3.2.} By $L\subset {\cal L}_{Lie}$ we denote the closed subspace of elements
$g=\sum_{n=1}^{\infty}g_{n}X_{n}t^{n}$, $g_{n}\in\Co$, $n\in\N$. Then there is a continuous linear isomorphism $A:{\cal L}_{S}\to L$ determined by $A(DL^{n-1}t^{n}):=X_{n}t^{n}$, $n\in\N$. Thus, $\Psi\circ A=id$. The map
$\Pi:=A\circ\Psi:{\cal L}_{Lie}\to L$ is a continuous linear projection onto $L$. Moreover, $id-\Pi: {\cal L}_{Lie}\to {\cal L}_{\widehat{\cal C}_{f}}$ is a 
continuous linear projection onto ${\cal L}_{\widehat{\cal C}_{f}}$. Hence
$\Pi\oplus (id-\Pi): {\cal L}_{Lie}\to L\oplus {\cal L}_{\widehat{\cal C}_{f}}$ is an isomorphism. Also, every element $g\in {\cal L}_{\widehat{\cal C}_{f}}$ is presented in the form
\begin{equation}\label{e3.23}
\begin{array}{c}
\displaystyle
g=\sum_{n=1}^{\infty}\left(\ \!\sum_{i_{1}+\cdots +i_{k}=n,\ k\geq 2}c_{i_{1},\dots, i_{k}}v_{i_{1},\dots, i_{k}}\right) t^{n}\ \ \ {\rm where}\\
\\
\displaystyle
v_{i_{1},\dots, i_{k}}:=[X_{i_{1}},[X_{i_{2}},[\ \cdots , [X_{i_{k-1}},X_{i_{k}}]\cdots\ \!]]]-\gamma_{i_{1},\dots, i_{k}}X_{n}.
\end{array}
\end{equation}
Recall that elements $\{v_{i_{1},\dots, i_{k}}\ :\ i_{1}+\cdots +i_{k}=n,\ k\geq 2\}$ are not linearly independent, cf. (\ref{e2.14}). The number of linearly independent elements in this set is 
$$
\frac{1}{n}\left(\ \!\sum_{d|n}(2^{n/d}-1)\cdot\mu(d)\right)-1.
$$
\begin{Proposition}\label{pr3.7}
$\widehat{\cal C}_{f}$ is the closure in $G_{f}(X)$ of the group $H$ generated by elements $\exp(c_{i_{1},\dots, i_{k}}v_{i_{1},\dots, i_{k}} t^{i_{1}+\cdots +i_{k}})$ for all possible
$v_{i_{1},\dots, i_{k}}$ and $c_{i_{1},\dots, i_{k}}\in\Co$.
\end{Proposition}
{\bf Proof.} According to Proposition \ref{pr2.22}, $\log(\overline{H})$ is the minimal closed Lie subalgebra of ${\cal L}_{Lie}$ containing all possible
elements
$v_{i_{1},\dots, i_{k}} t^{i_{1}+\cdots +i_{k}}$. By (\ref{e3.23}) this subalgebra coincides with
${\cal L}_{\widehat{\cal C}_{f}}$.
Thus, $\overline{H}=\widehat{\cal C}_{f}$.\ \ \ \ \ $\Box$\\
\\
{\bf Question 1.} {\em Is it true that every nontrivial element $\exp(c_{i_{1},\dots, i_{k}}v_{i_{1},\dots, i_{k}}t^{i_{1}+\cdots +i_{k}})$ belongs to $G_{f}(X)\setminus G(X)$, i.e., it cannot be presented by an element from $X$?}\\

{\bf 3.3.3.} Let us consider a continuous homomorphism $\pi_{2}:{\cal L}_{Lie}\to {\cal L}_{Lie}$ determined by the conditions
\begin{equation}\label{e3.24}
\pi_{2}(X_{i}):=\left\{
\begin{array}{ccc}
X_{i}&{\rm for}& i=1,2\\
\\
\frac{1}{i-2}\cdot [\pi_{2}(X_{i-1}), X_{1}]&{\rm for}&i\geq 3.
\end{array}
\right.
\end{equation}

From this definition we get
\begin{equation}\label{e3.25}
\pi_{2}(X_{i})=\frac{(-1)^{i}}{(i-2)!}\cdot\underbrace{[X_{1},[X_{1},[\ \!\cdots ,[X_{1},X_{2}]\cdots \ \!]]]}_{(i-2)-brackets}\ \ \ {\rm for}\ \ \ i\geq 3.
\end{equation}
Therefore $\pi_{2}$ maps ${\cal L}_{Lie}$ surjectively onto ${\cal L}_{Lie}^{2}$, the Lie algebra of the group $G_{f}(X^{2})$ of formal paths in $\Co^{2}$, see section 2.6.1.
\begin{Proposition}\label{pr3.8}
The kernel $Ker(\pi_{2})$ is a subgroup of ${\cal L}_{\widehat{\cal C}_{f}}$. It is the minimal normal closed Lie subalgebra ${\cal N}$ of ${\cal L}_{Lie}$ containing elements $\{s_{i} t^{i}\}_{i\geq 3}$ where
$$
s_{i}:=X_{i}-\frac{1}{i-2}\cdot [X_{i-1}, X_{1}],\ \ \ i\geq 3.
$$
\end{Proposition}
{\bf Proof.} By the definition each $s_{i}t^{i}\in Ker(\pi_{2})$. Since $\pi_{2}$ is continuous, ${\cal N}\subset Ker(\pi_{2})$.  Next, consider an element
$$
g=\sum_{n=1}^{\infty}\left(\ \!\sum_{i_{1}+\dots +i_{k}=n}c_{i_{1},\dots, i_{k}}[X_{i_{1}},[X_{i_{2}},[\ \cdots , [X_{i_{k-1}},X_{i_{k}}]\cdots \ ]]]\right) t^{n}=\sum_{n=1}^{\infty}g_{n}t^{n}\in {\cal L}_{Lie}.
$$
Using (\ref{e3.25}) one can represent each $g_{n}$ as $f_{n}+h_{n}$ where
$f_{n}t^{n}\in {\cal N}\cap (L_{n}\cdot t^{n})$ and $h_{n}t^{n}\in {\cal L}_{Lie}^{2}\cap (L_{n}\cdot t^{n})$
(here $L_{n}$ is the complex vector space generated by all possible $g_{n}$, see (\ref{e2.14})). Then we have
$$
g=f+h\ \ \ {\rm where}\ \ \ f:=\sum_{n=1}^{\infty}f_{n}t^{n}\in {\cal N}\ \ \ {\rm and}\ \ \ h:=\sum_{n=1}^{\infty}h_{n}t^{n}\in {\cal L}_{Lie}^{2}.
$$
In particular, $\pi_{2}(g)=h$. Thus $g\in Ker(\pi_{2})$ if and only if
$h=0$, i.e., $g\in {\cal N}$.

Further, each $s_{i}t^{i}\in {\cal L}_{\widehat{\cal C}_{f}}$, see (\ref{e3.21}). Hence
$Ker(\pi_{2})={\cal N}\subset {\cal L}_{\widehat{\cal C}_{f}}$.\ \ \ \ \ $\Box$

In the proof we have established also the natural decomposition 
\begin{equation}\label{e3.26}
{\cal L}_{Lie}={\cal N}\oplus {\cal L}_{Lie}^{2}.
\end{equation}
Taking the exponential map in (\ref{e3.26}) we obtain the following result.
\begin{Proposition}\label{pr3.9}
The group $G_{f}(X)$ is the semidirect product of the normal subgroup $\exp({\cal N})\subset\widehat{\cal C}_{f}$ and $G_{f}(X^{2})$. Moreover, $\exp({\cal N})$ is the  minimal normal closed subgroup of $G_{f}(X)$ containing elements $\exp(c_{i}s_{i}t^{i})$ for all possible $s_{i}$ and $c_{i}\in\Co$.
\end{Proposition}
{\bf Proof.} From the Campbell-Hausdorff formula for elements of $G_{f}(X)$ using the fact that ${\cal N}\subset {\cal L}_{Lie}$ is a closed ideal we obtain for $a\in {\cal N}$, $b\in {\cal L}_{Lie}^{2}$:
$$
e^{a+b}=e^{a}e^{b}e^{c_{1}}=e^{a}(e^{b}e^{c_{1}}e^{-b})e^{b}=(e^{a}e^{c_{2}})e^{b}=e^{c_{3}}e^{b}
$$
for some $c_{1},c_{2}, c_{3}\in {\cal N}$. This and (\ref{e3.26}) give the first statement of the proposition.

The second statement is the direct consequence of Proposition \ref{pr2.22} applied to elements of $R:=\{\log(e^{a}e^{s_{i}t^{i}}e^{-a})=e^{a}s_{i}t^{i}e^{-a}\in {\cal N}\ :\ a\in {\cal L}_{Lie},\ i\in\N\}$.\ \ \ \ \ $\Box$
\begin{R}\label{r3.10}
{\rm By the definition of $s_{i}$ each formal path $e^{c_{i}s_{i}t^{i}}$, $c_{i}\in\Co$, belongs to the {\em subgroup of formal paths in the subspace} $W_{i}\subset\Co^{\infty}$ where 
$$
W_{i}:=\{(z_{1},z_{2},\dots)\in\Co^{\infty}\ :\  z_{k}=0\ {\rm for\ all}\  k\not\in\{1,\ i-1,\ i\}\ \!\}\cong\Co^{3}.
$$
In particular, there exist elements $\gamma_{l}\in X$ with first integrals
$\widetilde\gamma_{l}: I_{T}\to W_{i}$, $l\in\N$, such that the sequence $\{\pi(\gamma_{l})\}_{l\in\N}\subset G(X)$ converges to $e^{c_{i}s_{i}t^{i}}$.
(Recall that we identify $G_{f}(X)$ with $\widehat E(G_{f}(X))$.)
}
\end{R}

Let 
$$
\Psi_{2}:=\Psi|_{G_{f}(X^{2})}: G_{f}(X^{2})\to S \ \!(\cong G[[r]])\ \ \
{\rm  and}\ \ \
\widehat P_{2}:=\widehat P|_{G_{f}(X^{2})}: G_{f}(X^{2})\to G[[r]].
$$
Then we have, cf. (\ref{e3.17}),
\begin{equation}\label{e3.27}
\widehat P_{2}=\Phi\circ\Psi_{2}.
\end{equation}

Next, we extend the homomorphism $\pi_{2}$ determined by (\ref{e3.24}) to a continuous endomorphism of the associative algebra ${\cal A}$, see (\ref{e2.6}). We retain the same symbol for the extension. Then $\pi_{2}$ maps $G_{f}(X)$ surjectively onto $G_{f}(X^{2})$.  Moreover, by Proposition \ref{pr3.9},  $Ker(\pi_{2}|_{G_{f}(X)})=\exp({\cal N})$. 
\begin{Proposition}\label{pr3.11}
The following identity is valid:
$$
\Psi(g)=(\Psi_{2}\circ\pi_{2})(g)\ \ \ {\rm for\ all}\ \ \ g\in G_{f}(X).
$$
\end{Proposition}
{\bf Proof.} Since for $h\in {\cal L}_{Lie}$ we have
$$
\Psi(e^{h})=e^{\Psi(h)}\ \ \ {\rm and}\ \ \ (\Psi_{2}\circ\pi_{2})(e^{h})=
e^{(\Psi_{2}\circ\pi_{2})(h)},
$$
it suffices to check the identity of the proposition for elements of ${\cal L}_{Lie}$.
Moreover, it suffices to check it for elements $X_{i}$, $i\in\N$ (because $\{X_{i}t^{i}\}_{i\in\N}$ are generators of ${\cal L}_{Lie}$). In this case we have by the definitions of  $\Psi$ and $\pi_{2}$ and by Lemma \ref{le3.2}
$$
\begin{array}{c}
\displaystyle
\Psi(X_{i})=DL^{i-1}\ \ \ {\rm and}\ \ \ (\Psi_{2}\circ\pi_{2})(X_{i})=DL^{i-1}\ \ \ {\rm for}\ \ \ i=1,2,\\
\\
\displaystyle
(\Psi_{2}\circ\pi_{2})(X_{i})=\frac{(-1)^{i}}{(i-2)!}\cdot\underbrace{ [D,[D,[\ \!\cdots ,[D,DL]\cdots \ \!]]]}_{(i-2)-brackets}=DL^{i-1}\ \ \ {\rm for}\ \ \ i\geq 3.
\end{array}
$$
This completes the proof of the proposition.\ \ \ \ \ $\Box$

From this proposition, (\ref{e3.17}) and (\ref{e3.27}) we obtain that
\begin{equation}\label{e3.28}
\widehat P=\widehat P_{2}\circ\pi_{2}.
\end{equation}
In particular, the homomorphism $\widehat P_{2}: G_{f}(X^{2})\to G[[r]]$ corresponding to the first return maps of ``generalized`` Abel equations is surjective. Moreover, Proposition \ref{pr3.9} implies that ${\cal C}_{f}$ is the semidirect product of $\exp({\cal N})$ and $\widehat{\cal C}_{f}^{\hspace*{0.5mm}2}:=\widehat{\cal C}_{f}\cap G_{f}(X^{2})$ (the {\em group of formal centers of Abel differential equations}).\\

{\bf 3.3.4.} It has been shown above that there is a reduction of the center problem for $G_{f}(X)$ to that for $G_{f}(X^{2})$. In this part we prove some results on the structure of the group $\widehat{\cal C}_{f}^{\hspace*{0.5mm}2}$.

First, observe that $\widehat{\cal C}_{f}^{\hspace*{0.5mm}2}\subset G_{f}(X^{2})$ is determined by systems of equations of form (\ref{e3.20}) in which all $i_{l}\in\{1,2\}$, $l\in\N$. In turn, the Lie algebra ${\cal L}_{\widehat{\cal C}_{f}^{\hspace*{0.5mm}2}}$ of
$\widehat{\cal C}_{f}^{\hspace*{0.5mm}2}$ is determined by the system equations of form (\ref{e3.21}) in which also all $i_{l}\in\{1,2\}$, $l\in\N$.

By $L_{2}\subset {\cal L}_{Lie}^{2}$ we denote the closed subspace of elements
$g=\sum_{n=1}^{\infty}g_{n}r_{n}t^{n}$, $g_{n}\in\Co$, $n\in\N$, where
$r_{1}=X_{1}$, $r_{2}=X_{2}$ and 
\begin{equation}\label{e3.29}
r_{n}:=\frac{(-1)^{n}}{(n-2)!}\cdot\underbrace{[X_{1},[X_{1},[\ \!\cdots ,[X_{1},X_{2}]\cdots \ \!]]]}_{(n-2)-brackets}\ \ \ {\rm for}\ \ \ n\geq 3.
\end{equation}
By the definition, $\Psi_{2}(r_{n})=DL^{n-1}$. Then there is a continuous linear isomorphism $A_{2}: {\cal L}_{S}\to L_{2}$ determined by
$A_{2}(DL^{n-1}t^{n}):=r_{n}t^{n}$, $n\in\N$, such that $\Psi_{2}\circ A_{2}=id$. The map $\Pi_{2}:=A_{2}\circ\Psi_{2}: {\cal L}_{Lie}^{2}\to L_{2}$ is a continuous linear projection onto $L_{2}$. Moreover, $id-\Pi_{2}: {\cal L}_{Lie}^{2}\to {\cal L}_{\widehat{\cal C}_{f}^{\hspace*{0.5mm}2}}$ is a continuous linear projection onto ${\cal L}_{\widehat{\cal C}_{f}^{\hspace*{0.5mm}2}}$. Hence, 
\begin{equation}\label{e3.30}
\Pi_{2}\oplus(id-\Pi_{2}):{\cal L}_{Lie}^{2}\to L_{2}\oplus {\cal L}_{\widehat{\cal C}_{f}^{\hspace*{0.5mm}2}}\ \ \  {\rm is\ an\ isomorphism}.
\end{equation}
This implies that every element
$g\in {\cal L}_{\widehat{\cal C}_{f}^{\hspace*{0.5mm}2}}$ is presented in the form
\begin{equation}\label{e3.31}
\begin{array}{c}
\displaystyle
g=\sum_{n=1}^{\infty}\left(\ \!\sum_{i_{1}+\cdots +i_{k}=n,\ n\geq 5}c_{i_{1},\dots, i_{k}}l_{i_{1},\dots, i_{k}}\right) t^{n}\ \ \ {\rm where\ all}\ \ \ i_{s}\in\{1,2\},\ s\in\N,\\
\\
\displaystyle
{\rm and}\ \ \ l_{i_{1},\dots, i_{k}}:=[X_{i_{1}},[X_{i_{2}},[\ \cdots , [X_{i_{k-1}},X_{i_{k}}]\cdots\ \!]]]-\gamma_{i_{1},\dots, i_{k}}r_{n}.
\end{array}
\end{equation}
The elements $\{l_{i_{1},\dots, i_{k}}\ :\ i_{1}+\cdots +i_{k}=n,\ n\geq 5,\ 
i_{s}\in\{1,2\}\}$ are not linearly independent. It follows from [M-KO, Theorem 3.1] that the number of linearly independent elements in this set is
\begin{equation}\label{e3.32}
\frac{1}{n}\left(\sum_{d|n}(\lambda_{1}^{n/d}+\lambda_{2}^{n/d})\cdot\mu(d)\right)-1
\end{equation}
where $\lambda_{1}=\frac{1+\sqrt{5}}{2}$, $\lambda_{2}=\frac{1-\sqrt{5}}{2}$, $\mu$ is the M\"{o}bius function, see section 2.4.1,
and the sum ranges over all integers which divide $n$.

Further, similarly to Propositions \ref{pr3.7} from (\ref{e3.30}) and (\ref{e3.31}) we get
\begin{Proposition}\label{pr3.12}
\begin{itemize}
\item[(1)]
There is a continuous map $T_{f}^{2}: G[[r]]\to G_{f}(X^{2})$ such that
$\widehat P_{2}\circ T_{f}^{2}=id$. Moreover, the map $\widetilde T_{f}^{2}: G[[r]]\times \widehat{\cal C}_{f}^{\hspace*{0.5mm}2}\to G_{f}(X^{2})$ defined by $\widetilde T_{f}^{2}(s,g)=T_{f}^{2}(s)\cdot g$ is a homeomorphism.
\item[(2)]
The group $\widehat{\cal C}_{f}^{\hspace*{0.5mm}2}$ is the closure in $G_{f}(X^{2})$ of the group $H_{2}$ generated by elements $\exp(c_{i_{1},\dots, i_{k}}l_{i_{1},\dots, i_{k}}t^{i_{1}+\cdots +i_{k}})$ for all possible $l_{i_{1},\dots, i_{k}}$ and
$c_{i_{1},\dots, i_{k}}\in\Co$ with $i_{1},\dots, i_{k}\in\{1,2\}$.
\end{itemize}
\end{Proposition}
{\bf Proof.} (1) We define the map $T_{f}^{2}$ by the formula
\begin{equation}\label{e3.33}
T_{f}^{2}:=\exp\circ A_{2}\circ\log\circ\Phi^{-1}.
\end{equation}
Then $T_{f}^{2}: G[[r]]\to G_{f}(X^{2})$ is continuous as the composite of continuous maps. Also, from (\ref{e3.27}) by the properties of $A_{2}$ we get
$$
\begin{array}{c}
\displaystyle
\widehat P_{2}\circ T_{f}^{2}=\Phi\circ\Psi_{2}\circ\exp\circ A_{2}\circ\log\circ\Phi^{-1}=\Phi\circ\exp\circ\Psi_{2}\circ A_{2}\circ\log\circ\Phi^{-1}=\\
\\
\displaystyle
\Phi\circ\exp\circ\ \! id\circ\log\circ\Phi^{-1}=\Phi\circ\Phi^{-1}=id.
\end{array}
$$

Now, for the map $\widetilde T_{f}^{2}: G[[r]]\times \widehat{\cal C}_{f}^{\hspace*{0.5mm}2}\to G_{f}(X^{2})$, $\widetilde T_{f}^{2}(s,g)=T_{f}^{2}(s)\cdot g$, we define the map $Q: G_{f}(X^{2})\to 
G[[r]]\times \widehat{\cal C}_{f}^{\hspace*{0.5mm}2}$ by the formula
\begin{equation}\label{e3.34}
Q(h):=(\widehat P_{2}(h),[(T_{f}^{2}\circ\widehat P_{2})(h)]^{-1}\cdot h).
\end{equation}
The second term here belongs to $\widehat{\cal C}_{f}^{\hspace*{0.5mm}2}$ because 
$$
\widehat P_{2}[(T_{f}^{2}(\widehat P_{2}(h)))^{-1}\cdot h]=
[(\widehat P_{2}\circ T_{f}^{2}\circ\widehat P_{2})(h)]^{-1}\cdot \widehat P_{2}(h)=[\widehat P_{2}(h)]^{-1}\cdot\widehat P_{2}(h)=1.
$$
Clearly, both $\widetilde T_{f}^{2}$ and $Q$ are continuous maps. Moreover,
$$
\begin{array}{c}
\displaystyle
(Q\circ\widetilde T_{f}^{2})(s,g)=Q(T_{f}^{2}(s)\cdot g)=(\widehat P_{2}(T_{f}^{2}(s)\cdot g), [(T_{f}^{2}\circ\widehat P_{2})(T_{f}^{2}(s)\cdot g)]^{-1}\cdot T_{f}^{2}(s)\cdot g)=\\
\\
\displaystyle
(s, [T_{f}^{2}(s)]^{-1}\cdot T_{f}^{2}(s)\cdot g)=(s,g).
\end{array}
$$
Thus $Q$ is the inverse map to $\widetilde T_{f}^{2}$, i.e., $\widetilde T_{f}^{2}$ is a homeomorphism.

(2) This statement is proved similarly to the proof of  Proposition \ref{pr3.7}.\ \ \ \ \ $\Box$\\
\\
{\bf Question 2.} {\em Is it true that there is a continuous map $T^{2}: G_{c}[[r]]\to G(X^{2})$ such that $\widehat P_{2}\circ T^{2}=id$?}\\

The affirmative answer to this question will show that each locally convergent series from $G_{c}[[r]]$ can be obtained as the first return map of an Abel differential equation. Moreover, as in the proof of Proposition \ref{pr3.6} we will get that the group of centers $\widehat{\cal C}^{\hspace*{0.5mm}2}:=\widehat{\cal C}_{f}^{\hspace*{0.5mm}2}\cap G(X^{2})$ of Abel differential equations is dense in the group of formal centers
$\widehat{\cal C}_{f}^{\hspace*{0.5mm}2}$. Also, $G(X^{2})$ will be homeomorphic to $G_{c}[[r]]\times \widehat{\cal C}^{\hspace*{0.5mm}2}$.\\

{\bf 3.3.5.} In this section we briefly discuss the center problem over a field $\F\subset\Co$.

Let $G_{\F}[[r]]\subset G[[r]]$ be the subgroup of formal power series with coefficients from $\F$. Let $I_{\F}^{k}$ be the normal subgroup of  $G_{\F}[[r]]$ consisting of elements $f$ of the form $f(z):=z+d_{k+1}z^{k+1}+d_{k+2}z^{k+2}+\cdots$. We equip $G_{\F}[[r]]$ with
$\{I_{\F}^{k}\}_{k\in\N}$-adic topology, i.e., a sequence $\{f_{i}\}_{i\in\N}\subset G_{\F}[[r]]$ converges to $f\in G_{\F}[[r]]$ if and only if for any $k\in\N$ there is a number $N_{k}\in\N$ such that for all $n\geq N_{k}$ the images of $f_{n}$ and $f$ in the quotient group $G_{\F}[[r]]/I_{\F}^{k}$ coincide. By $G_{c,\F}[[r]]\subset G_{\F}[[r]]$ we denote the subgroup of locally convergent power series in $G_{\F}[[r]]$ equipped with the induced topology.

Further, consider the groups $G(X_{\F})\subset G_{f}(X_{\F})$ of paths and formal paths over $\F$, see section 2.6.2. According to (\ref{e3.16}) the homomorphism $\widehat P$ (corresponding to the first return maps of equations (\ref{e1})) maps $G(X_{\F})$ and $G_{f}(X_{\F})$ into
$G_{c,\F}[[r]]$ and $G_{\F}[[r]]$, respectively. Also, by the definition of topologies on $G_{f}(X_{\F})$ and $G_{\F}[[r]]$, $\widehat P:G_{f}(X_{\F})\to G_{\F}[[r]]$ is a continuous homomorphism of topological groups.
We set $\widehat{\cal C}_{\F}:=G(X_{\F})\cap\widehat{\cal C}$ and 
$(\widehat{\cal C}_{\F})_{f}:=G_{f}(X_{\F})\cap\widehat{\cal C}_{f}$. These groups are referred to as the {\em groups of centers and formal centers over} $\F$. Then similarly to the results of the previous sections one can prove the following statements.
\begin{itemize}
\item[(1)]
$g\in (\widehat{\cal C}_{\F})_{f}$ if and only if the element $g\in G_{f}(X_{\F})$ satisfies equations (\ref{e3.20}).
\item[(2)]
The Lie algebra ${\cal L}_{(\widehat{\cal C}_{\F})_{f}}\subset {\cal L}_{Lie(\F)}$ of $(\widehat{\cal C}_{\F})_{f}$ consists of elements 
$$
\sum_{n=1}^{\infty}\left(\ \!\sum_{i_{1}+\dots +i_{k}=n}c_{i_{1},\dots, i_{k}}[X_{i_{1}},[X_{i_{2}},[\ \cdots , [X_{i_{k-1}},X_{i_{k}}]\cdots \ ]]]\right) t^{n}  
$$
with all $c_{i_{1},\dots, i_{k}}\in\F$ satisfying equations (\ref{e3.21}).
In particular the exponential map $\exp: {\cal L}_{(\widehat{\cal C}_{\F})_{f}}\to  (\widehat{\cal C}_{\F})_{f}$ is a homeomorphism.
\item[(3)]
\quad $\widehat{\cal C}_{\F}$ is a dense subgroup of $(\widehat{\cal C}_{\F})_{f}$.
\end{itemize}

The last statement is proved similarly to Proposition \ref{pr3.6} using Theorem 2.12 of [Br3]. This result asserts that there is a continuous embedding
$T_{\F}:G_{c,\F}[[r]]\to G(X_{\F})$ such that $\widehat P\circ T_{\F}=id$.
In particular, from here we obtain that $G(X_{\F})$ is homeomorphic to
$\widehat{\cal C}_{\F}\times G_{c,\F}[[r]]$ and $G_{f}(X_{\F})$ is homeomorphic to $(\widehat{\cal C}_{\F})_{f}\times G_{\F}[[r]]$.

Further, one can prove a version of Proposition \ref{pr3.7} (see Remark 2.23):
\begin{itemize}
\item[(4)]
$(\widehat{\cal C}_{\F})_{f}$ is the closure in $G_{f}(X_{\F})$ of the group $H_{\F}$ generated by elements \penalty-10000
$\exp(c_{i_{1},\dots, i_{k}}v_{i_{1},\dots, i_{k}}t^{i_{1}+\cdots +i_{k}})$ for all possible $v_{i_{1},\dots, i_{k}}$ determined by (\ref{e3.23}) and $c_{i_{1},\dots, i_{k}}\in\F$.
\end{itemize}

For the group ${\cal N}:=Ker(\pi_{2})$, see section 3.3.3, we set 
${\cal N}_{\F}:={\cal N}\cap {\cal L}_{Lie(\F)}$. Then repeating word-for-word the proof of Proposition \ref{pr3.8} we obtain that
\begin{itemize}
\item[(5)]
${\cal N}_{\F}$ is the minimal normal closed Lie subalgebra of ${\cal L}_{Lie(\F)}$ containing elements $s_{i}t^{i}=(X_{i}-\frac{1}{i-2}\cdot[X_{i-1},X_{1}])\ \!t^{i}$, $i\geq 3$. Moreover,
${\cal N}_{\F}\subset {\cal L}_{(\widehat{\cal C}_{\F})_{f}}$ and 
$$
{\cal L}_{Lie(\F)}={\cal N}_{\F}\oplus {\cal L}_{Lie(\F)}^{2}
$$
where ${\cal L}_{Lie(\F)}^{2}$ is the Lie algebra of $G_{f}(X_{\F}^{2}):=G_{f}(X_{\F})\cap G_{f}(X^{2})$.
\end{itemize}
From here as in the proof of Proposition \ref{pr3.9} we obtain
\begin{itemize}
\item[(6)]
The group $G_{f}(X_{\F})$ is the semidirect product of the normal subgroup $\exp({\cal N}_{\F})\subset (\widehat{\cal C}_{\F})_{f}$ and $G_{f}(X_{\F}^{2})$. Moreover, $\exp({\cal N}_{\F})$ is the minimal closed subgroup of $G_{f}(X_{\F})$ containing elements $\exp(c_{i}s_{i}t^{i})$ for all possible $s_{i}$ and $c_{i}\in\F$.
\end{itemize}

Finally, using the methods of the proof of Proposition \ref{pr3.12} one can prove similar statements for the group $(\widehat{\cal C}_{\F}^{2})_{f}:=(\widehat{\cal C}_{\F})_{f}\cap G_{f}(X_{\F}^{2})$:
\begin{itemize}
\item[(7)]
There is a continuous map $(T_{\F}^{2})_{f}:G_{\F}[[r]]\to G_{f}(X_{\F}^{2})$ such that $\widehat P_{2}\circ (T_{\F}^{2})_{f}=id$. 
The map $(\widetilde T_{\F}^{2})_{f}: G_{\F}[[r]]\times (\widehat{\cal C}_{\F}^{2})_{f}\to G_{f}(X_{\F}^{2})$ defined by $(\widetilde T_{\F}^{2})_{f}(s,g)=(T_{\F}^{2})_{f}(s)\cdot g$ is a homeomorphism. 
\item[(8)]
$(\widehat{\cal C}_{\F}^{2})_{f}$ is the closure in $G_{f}(X_{\F}^{2})$ of the group $(H_{\F})_{2}$ generated by elements $\exp(c_{i_{1},\dots, i_{k}}l_{i_{1},\dots, i_{k}}t^{i_{1}+\cdots +i_{k}})$ for all possible
$l_{i_{1},\dots, i_{k}}$ defined by (\ref{e3.31}) and $c_{i_{1},\dots, i_{k}}\in\F$,
$i_{1},\dots, i_{k}\in\{1,2\}$.
\end{itemize}

The main point of the proofs of (4)-(6) and (8) is that all elements $v_{i_{1},\dots, i_{k}}t^{i_{1}+\cdots +i_{k}}$, $s_{i}t^{i}$ and $l_{i_{1},\dots, i_{k}}t^{i_{1}+\cdots +i_{k}}$ belong to
${\cal L}_{Lie(\Q)}$ (and therefore to ${\cal L}_{Lie(\F)}$ for any field $\F\subset\Co$).
%========================
\subsect{\hspace*{-1em}. Group of Piecewise Linear Paths}

\quad Consider elements $g\in G_{f}(X)$ of the form 
\begin{equation}\label{e3.35}
g=e^{h}\ \ \ {\rm where}\ \ \ h=\sum_{i=1}^{\infty}c_{i}X_{i}t^{i},\ \ \
c_{i}\in\Co,\ i\in\N.
\end{equation}
By $PL\subset G_{f}(X)$ we denote the group generated by all such $g$. It will be called the {\em group of piecewise linear paths in} $\Co^{\infty}$.
\begin{R}\label{r3.13}
{\rm 
We can naturally extend the semigroup $X$ considering the set $\widetilde X:=(L^{\infty}(I_{T}))^{\N}$ of all possible sequences $a=(a_{1},a_{2},\dots)$, $a_{i}\in L^{\infty}(I_{T})$, $i\in\N$, with the multiplication defined in section 2.1.1. We consider each $L^{\infty}(I_{T})$ with the weak$*$ topology defined by $L^{1}(I_{T})$ (recall that $L^{\infty}(T)=(L^{1}(I_{T}))^{*}$) and
equip $\widetilde X$ with the corresponding product topology. 
Then $X$ is a dense subset of $\widetilde X$. Moreover, according to 
[Br3, Lemma 3.2] the quotient map $\pi: X\to G_{f}(X)$ is continuous and so is extended to a continuous map $\widetilde X\to G_{f}(X)$ (denoted also by $\pi$). Identifying $G_{f}(X)$ with its image under $\widehat{E}$, see (\ref{e2.21}), we obtain that $PL$ is the image under $\pi$ of the sub-semigroup $\widetilde X_{PL}$ of $\widetilde X$ generated by elements $c=(c_{1},c_{2},\dots)$, $c_{i}\in\Co$, $i\in\N$. In turn, first integrals of elements of this sub-semigroup are piecewise linear paths in $\Co^{\infty}$. This motivates the above definition.}
\end{R}

The group $PL$ contains the subgroup of rectangular paths $G(X_{rect})$, see section 2.4.2. In particular, $PL$ is a dense subgroup of $G_{f}(X)$, see also Proposition \ref{pr2.22}. One can also show (using, e.g., Theorem \ref{te2.2}) that $PL$ is isomorphic to the free $\Re$-product of groups $\Co$ (i.e., the set of generators of this product has the cardinality of the continuum). 

By $\widehat{\cal C}_{PL}:=\widehat {\cal C}_{f}\cap PL$ we denote the {\em group of formal centers in} $PL$. Then $\widehat{\cal C}_{PL}$ is the image in $G_{f}(X)$ of the semigroup ${\cal C}_{PL}\subset\widetilde X_{PL}$ consisting of all elements $a=(a_{1},a_{2},\dots)\in\widetilde X_{PL}$ such that monodromies of the equations
$$
H'(x)=\left(\ \!\sum_{i=1}^{\infty} a_{i}(x)DL^{i-1}t^{i}\right) H(x),\ \ \ x\in I_{T},
$$
are trivial.

Further, recall that there exist a continuous embedding $A:{\cal L}_{S}\to L$,
\penalty-10000
$L:=\{c\in {\cal L}_{Lie}\ :\ c=\sum_{i=1}^{\infty}c_{i}X_{i}t^{i},\ c_{i}\in\Co,\ i\in\N\}$, defined by $A(DL^{n-1}t^{n}):=X_{n}t^{n}$, so that $\Psi\circ A=id$, see section 3.3.2. For elements $a,b\in L$ we set
\begin{equation}\label{e3.36}
s(a,b):=(A\circ\Psi\circ\log)(e^{a}\cdot e^{b})\in L.
\end{equation}
Assume that
$$
a=\sum_{i=1}^{\infty}a_{i}X_{i}t^{i},\ \ \ b=\sum_{i=1}^{\infty}b_{i}X_{i}t^{i},\ \ \ a_{i}, b_{i}\in\Co, \ i\in\N.
$$
Using the Campbell-Hausdorff formula we have
$$
\log(e^{a}\cdot e^{b})=\sum_{n=1}^{\infty}\left(\ \!\sum_{i_{1}+\dots +i_{k}=n}s_{i_{1},\dots, i_{k}}(a,b)[X_{i_{1}},[X_{i_{2}},[\ \cdots , [X_{i_{k-1}},X_{i_{k}}]\cdots \ ]]]\right) t^{n}  
$$
where each $s_{i_{1},\dots, i_{k}}$ is a universal polynomial with rational coefficients in variables $a_{i_{1}},\dots ,a_{i_{k}}, b_{i_{1}},\dots, b_{i_{k}}$ such that $s_{i_{1},\dots, i_{k}}(a_{i_{1}}^{i_{1}},\dots, a_{i_{k}}^{i_{k}},b_{i_{1}}^{i_{1}},\dots, b_{i_{k}}^{i_{k}})$ is a homogeneous polynomial of degree $i_{1}+\cdots + i_{k}$. Then from (\ref{e3.23}), (\ref{e3.21}) we obtain
\begin{equation}\label{e3.37}
s(a,b)=\sum_{n=1}^{\infty}\left(\ \!\sum_{i_{1}+\cdots +i_{k}=n}\gamma_{i_{1},\dots, i_{k}}\cdot s_{i_{1},\dots, i_{k}}(a,b)\right)X_{n}t^{n}.
\end{equation}
In general, the complex vector space spanned by $a, b$ and $s(a,b)$ is $3$-dimensional.

Next, formula (\ref{e3.17}) implies that $\widehat P(e^{s(a,b)})=\widehat P(e^{a})\circ\widehat P(e^{b})$. In particular, \penalty-10000 $e^{a}\cdot e^{b}\cdot e^{-s(a,b)}\in\widehat{\cal C}_{PL}$.

Let us define a continuous embedding $T_{PL}: G[[r]]\to PL$ by the formula
$$
T_{PL}:=\exp\circ A\circ\log\circ\Phi^{-1},
$$
cf. (\ref{e3.33}). Then, $\widehat P\circ T_{PL}=id$.
\begin{Proposition}\label{pr3.14}
\begin{itemize}
\item[(1)]
The map $\widetilde T_{PL}:G[[r]]\times\widehat{\cal C}_{PL}\to PL$ defined by
$\widetilde T_{PL}(s,g)=T_{PL}(s)\cdot g$ is a homeomorphism.
\item[(2)]
$\widehat{\cal C}_{PL}$ is generated by elements $e^{a}\cdot e^{b}\cdot e^{-s(a,b)}$ for all possible $a,b\in L$.
\item[(3)]
$\widehat{\cal C}_{PL}$ is a dense subgroup of $\widehat{\cal C}_{f}$.
\end{itemize}
\end{Proposition}
{\bf Proof.} (1) The proof repeats literally the proof of Proposition \ref{pr3.12} (1).

(2) It follows easily from the definition of $T_{PL}$ and $\widehat P$, see (\ref{e3.17}), that
\begin{equation}\label{e3.38}
T_{PL}(\widehat P(e^{a}\cdot e^{b}))=e^{s(a,b)},\ \ \ a, b\in L.
\end{equation}
Suppose that $g=g_{1}\cdots g_{n}\in\widehat{\cal C}_{PL}$ where $g_{i}=e^{a_{i}}$, $a_{i}\in L$, $1\leq i\leq n$. We set 
$$
f_{i}:=\widehat P(g_{1}\cdots g_{i}),\ \ \ 1\leq i\leq n,
$$
where by the definition $f_{n}=1$. We also set for brevity
$$
c(h_{1},h_{2}):=h_{1}\cdot h_{2}\cdot e^{-s(\log(h_{1}),\log(h_{2}))},\ \ \ h_{1}, h_{2}\in PL.
$$
Then using (\ref{e3.38}) we obtain
$$
g=c(T_{PL}(f_{1}),g_{2})\cdot c(T_{PL}(f_{2}),g_{3})\cdots c(T_{PL}(f_{n-1}), g_{n})\cdot T_{PL}(f_{n}).
$$
Since $T_{PL}(f_{n})=1$, this implies the required statement.

(3) The proof repeats word-for-word the second part of the proof of Proposition \ref{pr3.6} and is based on the fact that $PL$ is dense in $G_{f}(X)$.\ \ \ \ \ $\Box$\\
\\
{\bf Question 3.} {\em Are there nontrivial elements in the group
$\widehat{\cal C}_{PL}^{n}:=\widehat{\cal C}_{f}^{n}\cap PL$ of piecewise linear centers in $\Co^{n}$? Here $\widehat{\cal C}_{f}^{n}:=\widehat{\cal C}_{f}\cap G_{f}(X^{n})$.}\\

We will return to this question in a forthcoming paper.
%========================

\end{document}